\documentclass[manuscript,printscheme,reqno]{aomart}
\usepackage{mystyle}
\usepackage{mathtools}
\usepackage{stmaryrd}
\usepackage{amssymb}
\usepackage{letltxmacro}

\newtheorem{thm}{Theorem}[section]
\newtheorem{cor}[thm]{Corollary}
\newtheorem{lem}[thm]{Lemma}
\newtheorem{prop}[thm]{Proposition}
\theoremstyle{definition}
\newtheorem{defn}[thm]{Definition}
\newtheorem{rem}[thm]{Remark}
\LetLtxMacro\origproof\proof
\LetLtxMacro\endorigproof\endproof
\makeatletter
\renewenvironment{proof}[1][\proofname]{\origproof[#1]\mbox{}\par\nobreak\@afterindentfalse\@afterheading\ignorespaces}{\endorigproof}
\makeatother

\DeclarePairedDelimiter{\floor}{\lfloor}{\rfloor}
\DeclarePairedDelimiter{\abs}{\lvert}{\rvert}
\DeclareMathOperator{\lcm}{lcm}
\newcommand{\N}{\mathbb{N}}
\newcommand{\Z}{\mathbb{Z}}
\newcommand{\R}{\mathbb{R}}
\newcommand{\bbinom}[3][b]{\genfrac{\llbracket}{\rrbracket}{0pt}{}{#2}{#3}_{\!#1}}

\allowdisplaybreaks

\title[Denumerant I]{Arithmetic Upper and Lower Bounds for the Denumerant Function I}
\author{Gerardo Ballesio}
\email{g.ballesio@gmail.com}
\givenname{Gerardo}
\surname{Ballesio}
\copyrightnote{}
\subject{primary}{matsc2020}{11P81}
\subject{secondary}{matsc2020}{05A17}
\subject{secondary}{matsc2020}{11B34}

\begin{document}

\begin{abstract}
	We use an elementary arithmetic argument, going back to Blom and Fr\"oberg, to obtain upper
	and lower bounds for Sylvester's denumerant function for arbitrary positive integers
	$a_1,\dotsc,a_k$. The bounds hold along the multiples of $\gcd(a_1,\dotsc,a_k)$, share the
	leading term of the denumerant, and yield its asymptotic behaviour;
	for large $n$ and under a mild condition, a polynomial refinement of the lower bound improves on
	the lower bound of Agnarsson.
\end{abstract}
\maketitle
\tableofcontents

% =====================================================================
\section*{Notation}\label{s:notation}
% =====================================================================
Throughout the paper we use the following notation.
\begin{itemize}
	\item $\N=\{0,1,2,\dotsc\}$ and, for $t\in\N$, $\N_{\geq t}=\{t,t+1,t+2,\dotsc\}$;
	      $\Z$ and $\R$ are the integers and the reals, and $\R_{\geq0}$ ($\R_{>0}$) the
	      non-negative (positive) reals.
	\item For $x,y\in\Z$ with $x\neq0$, $x\mid y$ means that $y=cx$ for some $c\in\Z$; here $y$
	      may be negative or zero.
	\item $\gcd$ and $\lcm$ are the greatest common divisor and the least common multiple of
	      positive integers.
	\item $\floor{x}$ is the floor of $x\in\R$, that is the integer such that
	      $\floor{x}\leq x<\floor{x}+1$. In particular $\floor{x}<0$ for $x<0$, so a sum
	      $\sum_{q=0}^{\floor{x}}$ is empty whenever $x<0$.
	\item For $m\in\Z$ and $j\in\N$, $\binom{m}{j}\coloneqq\frac{m(m-1)\dotsm(m-j+1)}{j!}$.
	\item For functions $f,g$ on $\N$ with $g\geq0$, $f(n)=O\bigl(g(n)\bigr)$ means that there are
	      $C\in\R_{>0}$ and $n_0\in\N$ with $\abs{f(n)}\leq C\,g(n)$ for every $n\geq n_0$.
	\item Empty sums are $0$ and empty products are $1$.
	\item $a_1,a_2,a_3,\dotsc$ is a fixed sequence of positive integers. Only finitely many of
	      them are involved in each statement; fixing the whole sequence once and for all simply
	      avoids carrying it around in the notation.
\end{itemize}

% =====================================================================
\section{Introduction}\label{s:intro}
% =====================================================================
In the nineteenth century, while studying the partition function,
J.~J.~Sylvester~\cite{sylvester:partition} and A.~Cayley~\cite{cayley} introduced
the \emph{denumerant}: the number of ways of writing a non-negative integer $n$ as a
non-negative integer combination of given positive integers $a_1,\dotsc,a_k$.

\begin{defn}\label{dfn:denumerant}
	For $k\in\N_{\geq1}$ and $n\in\N$, the \emph{denumerant} $D_k(n)$ is the number of
	solutions $(x_1,\dotsc,x_k)\in\N^k$ of $a_1x_1+\dotsb+a_kx_k=n$. Moreover, we set
	$d_k\coloneqq\gcd(a_1,\dotsc,a_k)$; in particular $d_1=a_1$.
\end{defn}

Since $d_k$ divides each of $a_1,\dotsc,a_k$, it divides every value of
$a_1x_1+\dotsb+a_kx_k$; so $D_k(n)=0$ unless $n$ is a multiple of $d_k$. This paper gives
explicit upper and lower bounds for $D_k(n)$.

\subsection{Known results}
We recall the known results in the notation of their authors; only $D_k(n)$ and $d_k$ are
needed to read them.

In 1962 G.~Blom and C.~Fr\"oberg~\cite{blom-froberg} proved, by elementary
methods, that if $a_1=1$ then for every $k\in\N_{\geq2}$ and every $n\in\N$
\begin{equation}\label{eq:intro-blom-froberg}
	\frac{n^{k-1}}{(k-1)!\,a_1\dotsm a_k}
	\;\leq\; D_k(n) \;\leq\;
	\frac{(n+s_k)^{k-1}}{(k-1)!\,a_1\dotsm a_k},
	\qquad
	s_k\coloneqq a_2+\frac12\sum_{i=3}^{k}a_i .
\end{equation}

In 2000 M.~B.~Nathanson~\cite{nathanson} showed, by arithmetic methods, that if
$a_1,\dotsc,a_k$ are distinct and $d_k=1$ then
\begin{equation}\label{eq:intro-nathanson}
	D_k(n)=\frac{n^{k-1}}{(k-1)!\,a_1\dotsm a_k}+O\bigl(n^{k-2}\bigr)
	\qquad(n\to+\infty).
\end{equation}

In 2002 G.~Agnarsson~\cite{agnarsson} proved, again by elementary methods, that
\begin{equation}\label{eq:intro-agnarsson}
	\frac{d_k}{a_1\dotsm a_k}\binom{n-B_k}{k-1}
	\;\leq\; D_k(n) \;\leq\;
	\frac{d_k}{a_1\dotsm a_k}\binom{n+A_k}{k-1},
\end{equation}
where $A_1=B_1=0$ and, for $k\geq2$,
\begin{equation}\label{eq:intro-AB}
	A_k=A_{k-1}+\frac{d_{k-1}}{d_k}\,a_k,
	\qquad
	B_k=B_{k-1}+\biggl(\frac{d_{k-1}}{d_k}-1\biggr)a_k-1.
\end{equation}
In~\cite{agnarsson} the upper bound is proved for every $n\in\N$ (Theorem~3.3), and the
lower bound for $n\geq B_k+k-1$ under the assumption $D_k(n)>0$ (Theorem~3.4),
which can be dropped when $d_k=1$ (Corollary~3.5).

A first version of the present work~\cite{ballesio-arxiv}, restricted to the case $d_k=1$ and
with slightly different constants, was posted on arXiv on 28 April 2022. That first
version already contains upper and lower bounds with the leading term
of~\eqref{eq:intro-nathanson}, a polynomial refinement of the lower bound and the asymptotic
formula~\eqref{eq:intro-nathanson} derived from them. Later, and independently,
K.~Gajdzica~\cite{gajdzica-jnt} proved, by comparison with integrals, that if
$a_1<a_2<\dotsb<a_k$, $\gcd(a_1,a_2)=1$ and $k\geq3$, then for every $n\in\N_{\geq1}$
\begin{equation}\label{eq:intro-gajdzica}
	\biggl|D_k(n)-\frac{n^{k-1}}{(k-1)!\,a_1\dotsm a_k}\biggr|<n^{k-2},
\end{equation}
and in the same paper (Lemma~5.3; see also~\cite{gajdzica-aa}) he obtained
two-sided bounds involving further terms of the
polynomial part of $D_k(n)$, assuming that, for a fixed $j\in\{2,\dotsc,k\}$, every multisubset
of $\{a_1,\dotsc,a_k\}$ with $j$ elements has greatest common divisor $1$.

Finally, recall that if $d_k=1$ every sufficiently large integer is representable,
and the largest non-representable one is the \emph{Frobenius number}
$g(a_1,\dotsc,a_k)$. In 1942 A.~Brauer~\cite{brauer} proved that
\begin{equation}\label{eq:intro-brauer}
	g(a_1,\dotsc,a_k)\;\leq\;\sum_{i=2}^{k}\frac{d_{i-1}}{d_i}\,a_i-\sum_{i=1}^{k}a_i .
\end{equation}

\subsection{Main result}
In~\eqref{eq:intro-AB} and~\eqref{eq:intro-brauer} the integers $a_i$ enter through the
products $\frac{d_{i-1}}{d_i}\,a_i$. Each of these products is a least common multiple,
namely $\lcm(d_{i-1},a_i)$, as we shall see in Lemma~\ref{lem:basic}~(ii). Our bounds are
written in terms of these least common multiples.

\begin{defn}\label{dfn:constants}
	For $k\in\N_{\geq2}$ we set
	\begin{equation*}
		b_k\coloneqq\lcm(d_{k-1},a_k),
		\qquad
		\Delta_k\coloneqq\sum_{i=2}^{k}(b_i-a_i),
		\qquad
		\nabla_k\coloneqq b_2+\frac12\sum_{i=3}^{k}b_i .
	\end{equation*}
\end{defn}

The number $\Delta_k$ is the shift of $n$ in our lower bound, and $\nabla_k$ is the shift of
$n$ in our upper bound. Since $D_k(n)=0$ unless $n$ is a multiple of $d_k$, we restrict to
such $n$. Our main result is the following.

\begin{thm}\label{thm:main}
	Let $k\in\N_{\geq2}$ and let $n\in\N$ be a multiple of $d_k$. Then
	\begin{equation}\label{eq:intro-main}
		\frac{d_k\,(n-\Delta_k)^{k-1}}{(k-1)!\,a_1\dotsm a_k}
		\;<\; D_k(n) \;\leq\;
		\frac{d_k\,(n+\nabla_k)^{k-1}}{(k-1)!\,a_1\dotsm a_k},
	\end{equation}
	where the lower bound holds for $n\geq\Delta_k$ and the upper bound for every such $n$.
	In particular, if $d_k=1$ then every $n\in\N$ is a multiple of $d_k$,
	and~\eqref{eq:intro-main} holds for every $n\in\N$ (the lower bound for $n\geq\Delta_k$).
\end{thm}

For $a_1=1$ the two bounds reduce to those of Blom and
Fr\"oberg~\eqref{eq:intro-blom-froberg}, and in general both have the leading term of
Nathanson's formula~\eqref{eq:intro-nathanson}.

\subsection{Organization of the paper}
Section~\ref{s:prelim} collects the preliminary results: the properties of $d_k$, $b_k$,
$\Delta_k$ and $\nabla_k$, a recursion for $D_k$, an estimate for sums of powers and a
binomial-like symbol built from elementary symmetric polynomials. Section~\ref{s:main} proves
Theorem~\ref{thm:main} together with a sharper, polynomial lower bound.
Section~\ref{s:asymptotic} derives the asymptotic formula and shows how the same technique,
simplified, can replace the argument of Nathanson. Section~\ref{s:comparison} compares our
bounds with the known results recalled above. Finally, Section~\ref{s:final} contains refined
bounds depending on the residue of $n$, and bounds for the number of solutions of
$a_1x_1+\dotsb+a_kx_k\leq n$.

% =====================================================================
\section{Preliminaries}\label{s:prelim}
% =====================================================================

\subsection{The integers \texorpdfstring{$d_k$, $b_k$, $\Delta_k$ and $\nabla_k$}{dk, bk, Delta k and nabla k}}
We collect the properties of the integers of Definitions~\ref{dfn:denumerant}
and~\ref{dfn:constants} that are used throughout the paper.

\begin{lem}\label{lem:basic}
	For every $k\in\N_{\geq2}$:
	\begin{enumerate}
		\item[\textup{(i)}] $d_k=\gcd(d_{k-1},a_k)$; in particular $d_k\mid d_{k-1}$;
		\item[\textup{(ii)}] $b_k=\frac{d_{k-1}}{d_k}\,a_k$; in particular $b_k$ is a
		      multiple of $a_k$ and $b_k\geq a_k$;
		\item[\textup{(iii)}] $b_2\dotsm b_k=\frac{a_1\dotsm a_k}{d_k}$;
		\item[\textup{(iv)}] $\Delta_k$ and $\nabla_k$ satisfy the recursion
		      \begin{equation}\label{eq:intro-recursion}
		      	\begin{aligned}
		      		\Delta_2&=b_2-a_2, &\qquad \Delta_{k+1}&=\Delta_k+b_{k+1}-a_{k+1},\\
		      		\nabla_2&=b_2,     &\qquad \nabla_{k+1}&=\nabla_k+\tfrac12\,b_{k+1},
		      	\end{aligned}
		      \end{equation}
		      and $\Delta_k$ is a non-negative integer;
		\item[\textup{(v)}] if $a_1=1$, then $d_k=1$ and $b_k=a_k$; hence $\Delta_k=0$;
		\item[\textup{(vi)}] $D_k(n)=0$ for every $n\in\N$ that is not a multiple of $d_k$
		      (this also holds for $k=1$).
	\end{enumerate}
\end{lem}
\begin{proof}
	\begin{enumerate}
		\item[\textup{(i)}] For $k=2$ both sides equal $\gcd(a_1,a_2)$, since $d_1=a_1$. For $k\geq3$,
		      an integer divides each of $a_1,\dotsc,a_{k-1}$ if and only if it divides
		      $d_{k-1}=\gcd(a_1,\dotsc,a_{k-1})$; hence $a_1,\dotsc,a_k$ and $d_{k-1},a_k$ have the
		      same common divisors, and so the same greatest one:
		      \begin{equation*}
		      	d_k=\gcd(a_1,\dotsc,a_k)=\gcd(d_{k-1},a_k).
		      \end{equation*}
		      Being a common divisor of $d_{k-1}$ and $a_k$, $d_k$ divides $d_{k-1}$.
		\item[\textup{(ii)}] Since $\gcd(x,y)\lcm(x,y)=xy$ for positive integers $x,y$, by part~(i)
		      \begin{equation*}
		      	b_k=\lcm(d_{k-1},a_k)=\frac{d_{k-1}a_k}{\gcd(d_{k-1},a_k)}=\frac{d_{k-1}}{d_k}\,a_k ,
		      \end{equation*}
		      and $\frac{d_{k-1}}{d_k}$ is a positive integer, again by part~(i).
		\item[\textup{(iii)}] By part~(ii) the product telescopes, and $d_1=a_1$:
		      \begin{equation*}
		      	b_2\dotsm b_k
		      	=\frac{d_1}{d_2}\,a_2\cdot\frac{d_2}{d_3}\,a_3\dotsm\frac{d_{k-1}}{d_k}\,a_k
		      	=\frac{d_1}{d_k}\,a_2\dotsm a_k
		      	=\frac{a_1a_2\dotsm a_k}{d_k}.
		      \end{equation*}
		\item[\textup{(iv)}] The sum defining $\Delta_2$ has the single term $i=2$, and the sum in
		      $\nabla_2$ is empty. Passing from $k$ to $k+1$, each of the two sums gains exactly its
		      term $i=k+1$, which gives the recursion. By part~(ii), each term $b_i-a_i$ of
		      $\Delta_k$ is a non-negative integer.
		\item[\textup{(v)}] If $a_1=1$ then $d_1=1$; since $d_k\mid d_{k-1}$ by part~(i), $d_k=1$ for
		      every $k$, and then $b_k=a_k$ by part~(ii). Hence every term of $\Delta_k$ vanishes.
		\item[\textup{(vi)}] If $(x_1,\dotsc,x_k)\in\N^k$ is a solution of $a_1x_1+\dotsb+a_kx_k=n$,
		      collecting $d_k$, which divides each of $a_1,\dotsc,a_k$, gives
		      \begin{equation*}
		      	n=d_k\biggl(\frac{a_1}{d_k}\,x_1+\dotsb+\frac{a_k}{d_k}\,x_k\biggr),
		      \end{equation*}
		      with an integer in brackets. Hence there are no solutions unless $d_k\mid n$.\qedhere
	\end{enumerate}
\end{proof}

\subsection{A recursion for the denumerant}

Passing from $k$ to $k+1$, we classify the solutions of $a_1x_1+\dotsb+a_{k+1}x_{k+1}=n$
according to the value $t=a_{k+1}x_{k+1}$, which ranges over the multiples of $a_{k+1}$ with
$0\leq t\leq n$. For a fixed value of $t$,
\begin{equation*}
	a_1x_1+\dotsb+a_kx_k+a_{k+1}x_{k+1}=n
	\iff
	a_1x_1+\dotsb+a_kx_k=n-a_{k+1}x_{k+1}=n-t,
\end{equation*}
and $x_{k+1}=\frac{t}{a_{k+1}}$ is determined by $t$; so the number of solutions with that
value of $t$ equals the number of solutions $(x_1,\dotsc,x_k)$ of $a_1x_1+\dotsb+a_kx_k=n-t$,
that is $D_k(n-t)$. Summing over $t$, and then dropping the terms in which $n-t$ is not a multiple
of $d_k$, which vanish by Lemma~\ref{lem:basic}~(vi), we get
\begin{equation}\label{eq:recursion-M}
	D_{k+1}(n)
	=\sum_{\substack{0\leq t\leq n\\ a_{k+1}\mid t}}D_k(n-t)
	=\sum_{\substack{0\leq t\leq n\\ a_{k+1}\mid t,\ d_k\mid n-t}}D_k(n-t).
\end{equation}
The next lemma shows that the values of $t$ in the last sum form an arithmetic progression with
difference $b_{k+1}$, and rewrites the sum accordingly.

\begin{lem}\label{lem:recursion}
	Let $k\in\N_{\geq1}$ and let $n\in\N$ be a multiple of $d_{k+1}$. Then there is exactly one
	integer $s$ such that
	\begin{equation*}
		0\leq s<b_{k+1},\qquad a_{k+1}\mid s,\qquad d_k\mid n-s ;
	\end{equation*}
	moreover $s\leq b_{k+1}-a_{k+1}$, and
	\begin{equation}\label{eq:recursion}
		D_{k+1}(n)=\sum_{q=0}^{\floor*{\frac{n-s}{b_{k+1}}}}D_k\bigl(n-s-qb_{k+1}\bigr),
	\end{equation}
	where every argument $n-s-qb_{k+1}$ is a non-negative multiple of $d_k$. In particular the
	sum is empty if $n<s$, and contains the term $q=0$ if $n\geq b_{k+1}-a_{k+1}$.
\end{lem}
\begin{proof}
	We prove, in this order: the existence of $s$; an equivalence that describes the values of
	$t$ in~\eqref{eq:recursion-M}; the uniqueness of $s$; the formula~\eqref{eq:recursion}.

	\begin{itemize}
		\item[$\scriptstyle\blacklozenge$] \emph{Step~1: existence of $s$.} Since $d_{k+1}=\gcd(d_k,a_{k+1})$ by
		Lemma~\ref{lem:basic}~(i), B\'ezout's identity gives $u,v\in\Z$ with
		$ua_{k+1}+vd_k=d_{k+1}$. Since $d_{k+1}\mid n$ by hypothesis, there is $c\in\N$ with
		$n=cd_{k+1}$; multiplying the previous identity by $c$,
		\begin{equation*}
			n=cd_{k+1}=cua_{k+1}+cvd_k .
		\end{equation*}
		Hence $t_0\coloneqq cua_{k+1}$ satisfies $a_{k+1}\mid t_0$ and $d_k\mid n-t_0=cvd_k$;
		$t_0$ may be negative, but only these two divisibilities will be used.
		Dividing $t_0$ by $b_{k+1}$ with remainder, we write $t_0=\ell b_{k+1}+s$ with $\ell\in\Z$ and
		$0\leq s<b_{k+1}$.
		Since $b_{k+1}=\lcm(d_k,a_{k+1})$ is a multiple of both $a_{k+1}$ and $d_k$,
		\begin{equation*}
			s=t_0-\ell b_{k+1}\ \text{is a multiple of } a_{k+1},
			\qquad
			n-s=(n-t_0)+\ell b_{k+1}\ \text{is a multiple of } d_k ,
		\end{equation*}
		so $s$ satisfies the three conditions of the statement. Moreover, $s$ and $b_{k+1}$ are both
		multiples of $a_{k+1}$ and $s<b_{k+1}$, hence $s\leq b_{k+1}-a_{k+1}$.

		\item[$\scriptstyle\blacklozenge$] \emph{Step~2: the values of $t$ form a progression.} Let $s$ be as in Step~1 and let $t$ be
		a multiple of $a_{k+1}$. We show that
		\begin{equation}\label{eq:recursion-equiv}
			d_k\mid n-t\iff b_{k+1}\mid t-s ,
		\end{equation}
		where $t-s$ may have any sign. We proceed in three parts.
		\begin{enumerate}
			\item[(a)] \emph{$d_k\mid n-t$ if and only if $d_k\mid t-s$.} Since $d_k\mid n-s$ by
			      Step~1, this follows from $n-t=(n-s)-(t-s)$: if $d_k$ divides two of the three
			      integers $n-t$, $n-s$, $t-s$, it divides the third.
			\item[(b)] \emph{$a_{k+1}\mid t-s$.} Indeed $a_{k+1}\mid t$ by assumption and
			      $a_{k+1}\mid s$ by Step~1. With~(a): $d_k\mid n-t$ if and only if $t-s$ is
			      divisible by both $a_{k+1}$ and $d_k$.
			\item[(c)] \emph{An integer $m$ is divisible by both $a_{k+1}$ and $d_k$ if and only if
			      $b_{k+1}\mid m$.} If $b_{k+1}\mid m$, then $a_{k+1}\mid m$ and $d_k\mid m$, since
			      $b_{k+1}$ is a multiple of both. Conversely, let $a_{k+1}\mid m$ and $d_k\mid m$, and
			      write $m=hb_{k+1}+\sigma$ with $h\in\Z$ and $0\leq\sigma<b_{k+1}$. Then
			      $\sigma=m-hb_{k+1}$ is divisible by both $a_{k+1}$ and $d_k$; if $\sigma$ were positive,
			      it would be a positive common multiple of $a_{k+1}$ and $d_k$ smaller than their
			      least common multiple $b_{k+1}$. Hence $\sigma=0$, that is $b_{k+1}\mid m$.
		\end{enumerate}
		Applying~(c) with $m=t-s$ and combining with~(b) gives~\eqref{eq:recursion-equiv}.

		\item[$\scriptstyle\blacklozenge$] \emph{Step~3: uniqueness of $s$.} Let $s'$ also satisfy the three conditions of the
		statement. Then $s'$ is a multiple of $a_{k+1}$ with $d_k\mid n-s'$,
		so~\eqref{eq:recursion-equiv} with $t=s'$ gives $b_{k+1}\mid s'-s$, that is,
		$b_{k+1}\mid\abs{s'-s}$. On the other hand $0\leq\abs{s'-s}<b_{k+1}$, because
		$0\leq s,s'<b_{k+1}$. The only multiple of $b_{k+1}$ in this range is $0$, so
		$\abs{s'-s}=0$, that is, $s'=s$.

		\item[$\scriptstyle\blacklozenge$] \emph{Step~4: the formula.} By~\eqref{eq:recursion-equiv}, the integers $t$ in the last sum
		of~\eqref{eq:recursion-M} are exactly those of the form $t=s+qb_{k+1}$ with $q\in\Z$ and
		$0\leq s+qb_{k+1}\leq n$; indeed each such $t$ is a multiple of $a_{k+1}$, since $a_{k+1}$
		divides both $s$ and $b_{k+1}$. Since $0\leq s<b_{k+1}$,
		\begin{equation*}
			s+qb_{k+1}\geq0\iff q\geq0,
			\qquad
			s+qb_{k+1}\leq n\iff q\leq\frac{n-s}{b_{k+1}}\iff q\leq\floor*{\frac{n-s}{b_{k+1}}}.
		\end{equation*}
		So $q$ ranges over the integers with $0\leq q\leq\floor*{\frac{n-s}{b_{k+1}}}$, and
		substituting $t=s+qb_{k+1}$ in~\eqref{eq:recursion-M} gives~\eqref{eq:recursion}. Each
		argument $n-t=n-s-qb_{k+1}$ is non-negative, since $t\leq n$, and a multiple of $d_k$, by
		Step~2. If $n<s$ no $q$ is admissible; if $n\geq b_{k+1}-a_{k+1}$, then $n\geq s$ by Step~1,
		and $q=0$ is admissible.
	\end{itemize}
\end{proof}

\subsection{Sums of powers}

We shall estimate sums of the form $\sum_{q=0}^{\floor{x}}(x-q)^j$ by telescoping: each term
$(x-q)^j$ is compared with an increment $F(x-q)-F(x-q-1)$ of a suitable function $F$, so that
the sum of the increments collapses to $F(x)$. The comparisons rest on the following pair of
identities, obtained by expanding $\bigl(u\pm\frac12\bigr)^p$ with the binomial theorem.
(They are the discrete counterparts of the midpoint and trapezoid comparisons between
$t^j$ and its integral, here obtained without integrals.)

\begin{lem}\label{lem:midpoint}
	Let $j\in\N_{\geq1}$. Then for every $u\in\R$
	\begin{align}
		\Bigl(u+\frac12\Bigr)^{j+1}-\Bigl(u-\frac12\Bigr)^{j+1}
		&=(j+1)\sum_{\substack{0\leq m\leq j\\ m\text{ even}}}\frac{1}{m+1}\binom{j}{m}\frac{u^{j-m}}{2^m},
		\label{eq:midpoint-diff}\\
		\Bigl(u+\frac12\Bigr)^{j}+\Bigl(u-\frac12\Bigr)^{j}
		&=2\sum_{\substack{0\leq m\leq j\\ m\text{ even}}}\binom{j}{m}\frac{u^{j-m}}{2^m}.
		\label{eq:midpoint-sum}
	\end{align}
	Consequently, for every $u\in\R_{\geq0}$
	\begin{equation}\label{eq:midpoint-ineq}
		(j+1)\,u^j
		\;\leq\;\Bigl(u+\frac12\Bigr)^{j+1}-\Bigl(u-\frac12\Bigr)^{j+1}
		\;\leq\;\frac{j+1}{2}\biggl(\Bigl(u+\frac12\Bigr)^{j}+\Bigl(u-\frac12\Bigr)^{j}\biggr).
	\end{equation}
\end{lem}
\begin{proof}
	By the binomial theorem,
	\begin{equation*}
		\Bigl(u+\frac12\Bigr)^{p}\pm\Bigl(u-\frac12\Bigr)^{p}
		=\sum_{m=0}^{p}\binom{p}{m}\frac{u^{p-m}}{2^m}\bigl(1\pm(-1)^m\bigr).
	\end{equation*}
	With the sign $+$ and $p=j$, only the even $m$ survive, each with factor $2$: this
	is~\eqref{eq:midpoint-sum}. With the sign $-$ and $p=j+1$, only the odd $m$ survive, each with
	factor $2$; writing $m=t+1$ with $t$ even, $0\leq t\leq j$, the term of index $m$ is
	\begin{equation*}
		2\binom{j+1}{t+1}\frac{u^{j-t}}{2^{t+1}}
		=\binom{j+1}{t+1}\frac{u^{j-t}}{2^{t}}
		=\frac{j+1}{t+1}\binom{j}{t}\frac{u^{j-t}}{2^{t}},
	\end{equation*}
	by the identity $\binom{j+1}{t+1}=\frac{j+1}{t+1}\binom{j}{t}$; this is~\eqref{eq:midpoint-diff}.

	Let now $u\geq0$, so that every term in~\eqref{eq:midpoint-diff} and~\eqref{eq:midpoint-sum}
	is non-negative. Keeping in~\eqref{eq:midpoint-diff} only the term $m=0$, which is
	$(j+1)u^j$, gives the first inequality in~\eqref{eq:midpoint-ineq}. Using
	$\frac{1}{m+1}\leq1$ in~\eqref{eq:midpoint-diff} and comparing with~\eqref{eq:midpoint-sum},
	\begin{equation*}
		\Bigl(u+\frac12\Bigr)^{j+1}-\Bigl(u-\frac12\Bigr)^{j+1}
		\leq(j+1)\sum_{\substack{0\leq m\leq j\\ m\text{ even}}}\binom{j}{m}\frac{u^{j-m}}{2^m}
		=\frac{j+1}{2}\biggl(\Bigl(u+\frac12\Bigr)^{j}+\Bigl(u-\frac12\Bigr)^{j}\biggr),
	\end{equation*}
	which is the second inequality.
\end{proof}

The first inequality in~\eqref{eq:midpoint-ineq} controls the terms $(x-q)^j$ with $x-q\geq1$
of a sum of powers; for the last term, which lies in $[0,1)$, we need the following inequality.

\begin{lem}\label{lem:last-term}
	Let $j\in\N_{\geq1}$. Then for every $t\in[0,1]$
	\begin{equation*}
		(j+1)\,t^j\leq\Bigl(t+\frac12\Bigr)^{j+1}.
	\end{equation*}
\end{lem}
\begin{proof}
	We prove the statement by induction on $j$.

	\emph{Base case $j=1$.} For $j=1$ the statement reads $2t\leq\bigl(t+\frac12\bigr)^2$, and
	\begin{equation*}
		\Bigl(t+\frac12\Bigr)^2-2t
		=t^2-t+\frac14
		=\Bigl(t-\frac12\Bigr)^2
		\geq0 .
	\end{equation*}

	\emph{Inductive step.} Assume the statement for some $j\geq1$, that is,
	$(j+1)\,t^j\leq\bigl(t+\frac12\bigr)^{j+1}$ for every $t\in[0,1]$. We have to prove that
	$(j+2)\,t^{j+1}\leq\bigl(t+\frac12\bigr)^{j+2}$ for every $t\in[0,1]$. Let $t\in[0,1]$.
	Multiplying the inductive hypothesis by $t+\frac12>0$ for the first inequality, and using
	$t\leq1\leq\frac{j+1}{2}$ and $t^j\geq0$ for the last one,
	\begin{equation*}
		\Bigl(t+\frac12\Bigr)^{j+2}
		=\Bigl(t+\frac12\Bigr)\Bigl(t+\frac12\Bigr)^{j+1}
		\geq\Bigl(t+\frac12\Bigr)(j+1)\,t^j
		=(j+2)\,t^{j+1}+\biggl(\frac{j+1}{2}-t\biggr)t^j
		\geq(j+2)\,t^{j+1},
	\end{equation*}
	which is the statement for $j+1$.
\end{proof}

We can now estimate the sums of powers.

\begin{lem}\label{lem:sum-bounds}
	For every $j\in\N_{\geq1}$ and every $x\in\R_{\geq0}$
	\begin{equation*}
		\frac{x^{j+1}}{j+1}+\frac{x^j}{2}
		\;\leq\;\sum_{q=0}^{\floor{x}}(x-q)^j
		\;\leq\;\frac{\bigl(x+\frac12\bigr)^{j+1}}{j+1}.
	\end{equation*}
\end{lem}
\begin{proof}
	Let $t_q\coloneqq x-q$ for $0\leq q\leq\floor{x}$, so that
	\begin{equation*}
		\sum_{q=0}^{\floor{x}}(x-q)^j=t_0^{\,j}+t_1^{\,j}+\dotsb+t_{\floor{x}}^{\,j}.
	\end{equation*}
	The first term is $t_0=x$, and each term is the previous one minus $1$:
	\begin{equation}\label{eq:t-step}
		t_{q+1}
		=x-(q+1)
		=(x-q)-1
		=t_q-1 .
	\end{equation}
	If $q<\floor{x}$, then $q\leq\floor{x}-1\leq x-1$, hence $t_q=x-q\geq1$; the last term,
	$t_{\floor{x}}=x-\floor{x}$, lies in $[0,1)$.

	\emph{The method.} The sum of the powers $t_q^{\,j}$ has no simple closed form, so we compare it
	with a sum that collapses. Let $F$ be a real function on $\R$. By~\eqref{eq:t-step}, for
	$q<\floor{x}$ the difference $F(t_q)-F(t_{q+1})=F(t_q)-F(t_q-1)$ is the increment of $F$ between
	$t_q-1$ and $t_q$, and the sum of these increments telescopes:
	\begin{equation}\label{eq:telescope}
		\sum_{q=0}^{\floor{x}-1}\bigl(F(t_q)-F(t_{q+1})\bigr)
		=\bigl(F(t_0)-F(t_1)\bigr)+\bigl(F(t_1)-F(t_2)\bigr)+\dotsb+\bigl(F(t_{\floor{x}-1})-F(t_{\floor{x}})\bigr)
		=F(t_0)-F(t_{\floor{x}}).
	\end{equation}
	Suppose now that $F$ satisfies the two conditions
	\begin{align}
		t_q^{\,j}&\leq F(t_q)-F(t_q-1)\quad\text{for every } q<\floor{x},\label{eq:method-a}\\
		t_{\floor{x}}^{\,j}&\leq F(t_{\floor{x}}).\label{eq:method-b}
	\end{align}
	Then, summing~\eqref{eq:method-a} over $q<\floor{x}$, using~\eqref{eq:telescope} and
	adding~\eqref{eq:method-b},
	\begin{equation}\label{eq:method-c}
		\sum_{q=0}^{\floor{x}}t_q^{\,j}
		\leq\bigl(F(t_0)-F(t_{\floor{x}})\bigr)+F(t_{\floor{x}})
		=F(t_0)
		=F(x),
	\end{equation}
	which depends only on $x$. In the same way, if~\eqref{eq:method-a} and~\eqref{eq:method-b} hold
	with $\geq$ in place of $\leq$, then~\eqref{eq:method-c} holds with $\geq$ in place of $\leq$.

	So far~\eqref{eq:method-a} and~\eqref{eq:method-b} are assumptions on $F$: nothing has been
	proved for a specific function. In the rest of the proof we consider the two functions
	\begin{equation*}
		L(t)\coloneqq\frac{t^{j+1}}{j+1}+\frac{t^j}{2},
		\qquad
		U(t)\coloneqq\frac{\bigl(t+\frac12\bigr)^{j+1}}{j+1},
	\end{equation*}
	for which $L(x)$ and $U(x)$ are the two bounds of the statement. We prove that $F=U$
	satisfies~\eqref{eq:method-a} and~\eqref{eq:method-b}, so that~\eqref{eq:method-c} gives the
	upper bound; and that $F=L$ satisfies them with $\geq$, so that~\eqref{eq:method-c} with $\geq$
	gives the lower bound.

	\emph{Upper bound.} For $q<\floor{x}$, the first inequality in~\eqref{eq:midpoint-ineq} with
	$u=t_q\geq0$, divided by $j+1$, gives
	\begin{equation*}
		t_q^{\,j}
		\leq\frac{\bigl(t_q+\frac12\bigr)^{j+1}-\bigl(t_q-\frac12\bigr)^{j+1}}{j+1}
		=U(t_q)-U(t_q-1).
	\end{equation*}
	This is condition~\eqref{eq:method-a} for $F=U$. For the last term, Lemma~\ref{lem:last-term} with $t=t_{\floor{x}}\in[0,1)$, divided by $j+1$,
	gives
	\begin{equation*}
		t_{\floor{x}}^{\,j}
		\leq\frac{\bigl(t_{\floor{x}}+\frac12\bigr)^{j+1}}{j+1}
		=U(t_{\floor{x}}).
	\end{equation*}
	This is condition~\eqref{eq:method-b} for $F=U$. Since $U$ satisfies both conditions,
	\eqref{eq:method-c} holds for $F=U$:
	\begin{equation*}
		\sum_{q=0}^{\floor{x}}t_q^{\,j}
		\leq U(x)
		=\frac{\bigl(x+\frac12\bigr)^{j+1}}{j+1}.
	\end{equation*}

	\emph{Lower bound.} For $q<\floor{x}$ we have $t_q\geq1$, so $u\coloneqq t_q-\frac12\geq\frac12\geq0$,
	and $u+\frac12=t_q$, $u-\frac12=t_q-1$. The second inequality in~\eqref{eq:midpoint-ineq} with
	this $u$, divided by $j+1$, gives
	\begin{equation*}
		\frac{t_q^{\,j+1}-(t_q-1)^{j+1}}{j+1}\leq\frac{t_q^{\,j}+(t_q-1)^{j}}{2}.
	\end{equation*}
	Adding $\frac{t_q^{\,j}-(t_q-1)^{j}}{2}$ to both sides, by the definition of $L$,
	\begin{equation*}
		L(t_q)-L(t_q-1)
		=\frac{t_q^{\,j+1}-(t_q-1)^{j+1}}{j+1}+\frac{t_q^{\,j}-(t_q-1)^{j}}{2}
		\leq\frac{t_q^{\,j}+(t_q-1)^{j}}{2}+\frac{t_q^{\,j}-(t_q-1)^{j}}{2}
		=t_q^{\,j}.
	\end{equation*}
	This is condition~\eqref{eq:method-a} with $\geq$ for $F=L$. For the last term, since
	$0\leq t_{\floor{x}}<1$ and $j\geq1$, we have
	$\frac{t_{\floor{x}}}{j+1}\leq\frac12$, hence
	\begin{equation*}
		L(t_{\floor{x}})
		=t_{\floor{x}}^{\,j}\biggl(\frac{t_{\floor{x}}}{j+1}+\frac12\biggr)
		\leq t_{\floor{x}}^{\,j}\biggl(\frac{1}{2}+\frac12\biggr)
		=t_{\floor{x}}^{\,j}.
	\end{equation*}
	This is condition~\eqref{eq:method-b} with $\geq$ for $F=L$. Since $L$ satisfies both
	conditions with $\geq$, \eqref{eq:method-c} holds with $\geq$ for $F=L$:
	\begin{equation*}
		\sum_{q=0}^{\floor{x}}t_q^{\,j}
		\geq L(x)
		=\frac{x^{j+1}}{j+1}+\frac{x^j}{2}.\qedhere
	\end{equation*}
\end{proof}

\begin{rem}\label{rem:last-term}
	In the proof of Lemma~\ref{lem:sum-bounds} the terms $t_q^{\,j}$ with $q<\floor{x}$ are compared
	with increments of $F$, while the last term $t_{\floor{x}}^{\,j}$ is compared with the single
	value $F(t_{\floor{x}})$. Suppose we used an increment also for the last term, that is, the
	inequality between $t_{\floor{x}}^{\,j}$ and $F(t_{\floor{x}})-F(t_{\floor{x}}-1)$. Then the
	telescoping sum~\eqref{eq:telescope} would have one more term, and the final result would be
	$F(t_0)-F(t_{\floor{x}}-1)$ instead of $F(t_0)$. This does not work, for a different reason in
	each of the two bounds.

	\emph{Upper bound, $F=U$.} We would get
	\begin{equation*}
		\sum_{q=0}^{\floor{x}}t_q^{\,j}\leq U(x)-U(t_{\floor{x}}-1),
		\qquad\text{where}\qquad
		U(t_{\floor{x}}-1)=\frac{\bigl(t_{\floor{x}}-\frac12\bigr)^{j+1}}{j+1}.
	\end{equation*}
	This gives the bound $U(x)$ of the lemma only if $U(t_{\floor{x}}-1)\geq0$. But if
	$t_{\floor{x}}<\frac12$ and $j$ is even, then $\bigl(t_{\floor{x}}-\frac12\bigr)^{j+1}$ is an
	odd power of a negative number, so $U(t_{\floor{x}}-1)<0$.

	\emph{Lower bound, $F=L$.} For $q<\floor{x}$ the inequality $t_q^{\,j}\geq L(t_q)-L(t_q-1)$ was
	obtained from Lemma~\ref{lem:midpoint} with $u=t_q-\frac12$, and that lemma requires $u\geq0$,
	that is $t_q\geq\frac12$. The last term $t_{\floor{x}}$ lies in $[0,1)$ and may be smaller than
	$\frac12$; in that case Lemma~\ref{lem:midpoint} does not apply.
\end{rem}

\begin{rem}\label{rem:sum-bounds-bf}
	Lemma~\ref{lem:sum-bounds} is the lemma in Section~2.2 of Blom and
	Fr\"oberg~\cite{blom-froberg}, who prove the lower bound together with the sharper upper bound
	$\frac{x^{j+1}}{j+1}+\frac{x^j}{2}+\frac{j\,x^{j-1}}{8}$, which consists of the first three terms of
	the binomial expansion of $\frac{(x+\frac12)^{j+1}}{j+1}$. They use the lower bound for their
	refined inequality~(B), and the upper bound only after replacing it by
	$\frac{(x+\frac12)^{j+1}}{j+1}$, in the proof of their inequality~(A); so the two lemmas are used
	in the same way. Their proof is by induction on $j$
	and uses differentiation with respect to $x$; the proof above uses only the binomial theorem
	and a telescoping sum.
\end{rem}

\subsection{Elementary symmetric polynomials}

For $m\in\N$ and $i\in\Z$, the \emph{elementary symmetric polynomial} of degree $i$ in
$X_1,\dotsc,X_m$ is
\begin{equation*}
	e_i(X_1,\dotsc,X_m)\coloneqq
	\begin{cases}
		\displaystyle\sum_{1\leq l_1<\dotsb<l_i\leq m}X_{l_1}\dotsm X_{l_i}, & 0\leq i\leq m,\\[2ex]
		0, & i<0 \text{ or } i>m.
	\end{cases}
\end{equation*}
For $i=0$ the sum has a single term, the empty product, so $e_0=1$, also for $m=0$.
To shorten the formulas, for $m\in\N$ and $i\in\Z$ we write
\begin{equation*}
	\bbinom{m}{i}\coloneqq e_i\Bigl(\frac{b_3}{2},\dotsc,\frac{b_{m+2}}{2}\Bigr).
\end{equation*}
The next lemma shows that this symbol behaves like a binomial coefficient.

\begin{lem}\label{lem:sym-poly}
	For every $m\in\N$ and $i\in\Z$:
	\begin{enumerate}
		\item[\textup{(i)}] $\bbinom{m}{0}=1$, $\bbinom{m}{i}=0$ for $i<0$ or $i>m$, and
		      $\bbinom{m}{i}\geq0$;
		\item[\textup{(ii)}] if $m\geq1$, then
		      $\bbinom{m}{i}=\bbinom{m-1}{i}+\frac{b_{m+2}}{2}\bbinom{m-1}{i-1}$;
		\item[\textup{(iii)}] if $b_3=\dotsb=b_{m+2}=2$ and $i\geq0$, then $\bbinom{m}{i}=\binom{m}{i}$;
		\item[\textup{(iv)}] more generally, let $c=(c_3,c_4,c_5,\dotsc)$ be a sequence of real
		      numbers and set
		      $\bbinom[c]{m}{i}\coloneqq e_i\bigl(\frac{c_3}{2},\dotsc,\frac{c_{m+2}}{2}\bigr)$. Then
		      $\bbinom[c]{m}{0}=1$, and $\bbinom[c]{m}{i}=0$ for $i<0$ or $i>m$; if $c_l\geq0$ for
		      every $l$, then $\bbinom[c]{m}{i}\geq0$; and if $m\geq1$, then
		      $\bbinom[c]{m}{i}=\bbinom[c]{m-1}{i}+\frac{c_{m+2}}{2}\bbinom[c]{m-1}{i-1}$.
	\end{enumerate}
\end{lem}
\begin{proof}
	Set $X_l\coloneqq\frac{b_{l+2}}{2}$ for $1\leq l\leq m$. Then
	\begin{equation*}
		\bbinom{m}{i}=e_i(X_1,\dotsc,X_m)
		\qquad\text{and, for } m\geq1,\qquad
		\bbinom{m-1}{i}=e_i(X_1,\dotsc,X_{m-1}).
	\end{equation*}
	\begin{enumerate}
		\item[\textup{(i)}] By definition $e_0=1$ and $e_i=0$ for $i<0$ or $i>m$; for $0\leq i\leq m$,
		      $e_i(X_1,\dotsc,X_m)$ is a sum of products of the positive numbers $X_l$.
		\item[\textup{(ii)}] We show that
		      \begin{equation*}
		      	e_i(X_1,\dotsc,X_m)=e_i(X_1,\dotsc,X_{m-1})+X_m\,e_{i-1}(X_1,\dotsc,X_{m-1}).
		      \end{equation*}
		      For $1\leq i\leq m$, split the products $X_{l_1}\dotsm X_{l_i}$ according to whether
		      they contain the factor $X_m$ or not: those that do not contain it give
		      $e_i(X_1,\dotsc,X_{m-1})$ (which is $0$ if $i=m$), and those that contain it give
		      $X_m\,e_{i-1}(X_1,\dotsc,X_{m-1})$. For $i=0$ both sides equal $1$, since $e_{-1}=0$;
		      for $i<0$ or $i>m$ both sides equal $0$. Since $X_m=\frac{b_{m+2}}{2}$, this is~(ii).
		\item[\textup{(iii)}] If every $X_l$ equals $1$, each product $X_{l_1}\dotsm X_{l_i}$ equals
		      $1$, so $e_i(1,\dotsc,1)$ is the number of subsets $\{l_1,\dotsc,l_i\}$ of
		      $\{1,\dotsc,m\}$, that is $\binom{m}{i}$ (which is $0$ for $i>m$).
		\item[\textup{(iv)}] The proofs of~(i) and~(ii) use the numbers $X_l$ only through the
		      definition of $e_i$, except for the inequality in~(i), which uses that they are
		      non-negative. So they hold with $X_l\coloneqq\frac{c_{l+2}}{2}$ in place of
		      $\frac{b_{l+2}}{2}$, the inequality under the assumption $c_l\geq0$.\qedhere
	\end{enumerate}
\end{proof}
Part~(iv) is not needed in this paper, where only the sequence $b_3,b_4,\dotsc$ is used; we state
it for later use.

% =====================================================================
\section{Main results}\label{s:main}
% =====================================================================
Throughout this section the coefficient of the bounds is written as $\frac{1}{b_2\dotsm b_k}$,
which equals $\frac{d_k}{a_1\dotsm a_k}$ by Lemma~\ref{lem:basic}~(iii).

\subsection{The case \texorpdfstring{$k=2$}{k=2}}

For $k=2$ the recursion of Lemma~\ref{lem:recursion} counts the terms of an arithmetic
progression, and both bounds follow at once.

\begin{prop}\label{prop:two}
	Let $n\in\N$ be a multiple of $d_2$, and let $s$ be as in Lemma~\ref{lem:recursion} with $k=1$.
	Then
	\begin{equation*}
		D_2(n)=
		\begin{cases}
			\floor*{\frac{n-s}{b_2}}+1, & n\geq s,\\[1ex]
			0, & n<s.
		\end{cases}
	\end{equation*}
	Consequently $D_2(n)\leq\frac{n+\nabla_2}{b_2}$, and $D_2(n)>\frac{n-\Delta_2}{b_2}$ if
	$n\geq\Delta_2$.
\end{prop}
\begin{proof}
	\emph{The value of $D_1$.} Let $m\in\N$ be a multiple of $d_1=a_1$. The equation $a_1x_1=m$ has
	exactly one solution $x_1\in\N$, namely $x_1=\frac{m}{a_1}$. Hence
	\begin{equation}\label{eq:D1}
		D_1(m)=1\quad\text{for every multiple } m\in\N \text{ of } d_1 .
	\end{equation}

	\emph{The formula for $D_2(n)$.} By Lemma~\ref{lem:recursion} with $k=1$,
	\begin{equation*}
		D_2(n)=\sum_{q=0}^{\floor*{\frac{n-s}{b_2}}}D_1(n-s-qb_2),
	\end{equation*}
	and every argument $n-s-qb_2$ is a non-negative multiple of $d_1$. By~\eqref{eq:D1} every term
	of the sum equals $1$, so $D_2(n)$ is the number of terms of the sum, that is, the number of
	integers $q$ with $0\leq q\leq\floor*{\frac{n-s}{b_2}}$. If $n<s$, then $\frac{n-s}{b_2}<0$, so
	$\floor*{\frac{n-s}{b_2}}<0$ and there is no such $q$: $D_2(n)=0$. If $n\geq s$, then
	$\floor*{\frac{n-s}{b_2}}\geq0$ and these integers are $0,1,\dotsc,\floor*{\frac{n-s}{b_2}}$:
	$D_2(n)=\floor*{\frac{n-s}{b_2}}+1$.

	\emph{Upper bound.} We have to prove that $D_2(n)\leq\frac{n+\nabla_2}{b_2}$. By
	Lemma~\ref{lem:basic}~(iv), $\nabla_2=b_2$. If $n<s$, then $D_2(n)=0$, and
	$\frac{n+b_2}{b_2}>0$ because $n\geq0$ and $b_2\geq1$. If $n\geq s$, then, using
	$\floor{y}\leq y$ for the first inequality and $s\geq0$ for the second one,
	\begin{equation*}
		D_2(n)
		=\floor*{\frac{n-s}{b_2}}+1
		\leq\frac{n-s}{b_2}+1
		\leq\frac{n}{b_2}+1
		=\frac{n+b_2}{b_2}
		=\frac{n+\nabla_2}{b_2}.
	\end{equation*}

	\emph{Lower bound.} Let $n\geq\Delta_2$. We have to prove that $D_2(n)>\frac{n-\Delta_2}{b_2}$.
	By Lemma~\ref{lem:basic}~(iv), $\Delta_2=b_2-a_2$, and by Lemma~\ref{lem:recursion} with $k=1$,
	$s\leq b_2-a_2$; hence $s\leq\Delta_2\leq n$, and $D_2(n)=\floor*{\frac{n-s}{b_2}}+1$. Using
	$\floor{y}+1>y$ for the first inequality and $s\leq\Delta_2$ for the second one,
	\begin{equation*}
		D_2(n)
		=\floor*{\frac{n-s}{b_2}}+1
		>\frac{n-s}{b_2}
		\geq\frac{n-\Delta_2}{b_2}.\qedhere
	\end{equation*}
\end{proof}

\subsection{Upper bound}

\begin{thm}\label{thm:upper}
	Let $k\in\N_{\geq2}$ and let $n\in\N$ be a multiple of $d_k$. Then
	\begin{equation*}
		D_k(n)\leq\frac{(n+\nabla_k)^{k-1}}{(k-1)!\,b_2\dotsm b_k}.
	\end{equation*}
\end{thm}
\begin{proof}
	We prove the statement by induction on $k$.

	\emph{Base case $k=2$.} For $k=2$ the statement reads $D_2(n)\leq\frac{n+\nabla_2}{b_2}$, which is
	the upper bound of Proposition~\ref{prop:two}.

	\emph{Inductive step.} Assume the statement for some $k\geq2$, that is, for every multiple
	$m\in\N$ of $d_k$,
	\begin{equation}\label{eq:upper-hyp}
		D_k(m)\leq\frac{(m+\nabla_k)^{k-1}}{(k-1)!\,b_2\dotsm b_k}.
	\end{equation}
	Let $n\in\N$ be a multiple of $d_{k+1}$, and let $s$ be the integer of Lemma~\ref{lem:recursion}.
	We have to prove that
	\begin{equation}\label{eq:upper-goal}
		D_{k+1}(n)\leq\frac{(n+\nabla_{k+1})^{k}}{k!\,b_2\dotsm b_{k+1}}.
	\end{equation}

	\emph{The case $n<s$.} By Lemma~\ref{lem:recursion} the sum in~\eqref{eq:recursion} is empty, so
	$D_{k+1}(n)=0$. The right-hand side of~\eqref{eq:upper-goal} is non-negative, because
	$n\geq0$ and $\nabla_{k+1}\geq0$ (by Definition~\ref{dfn:constants}, $\nabla_{k+1}$ is a sum of
	positive numbers). Hence~\eqref{eq:upper-goal} holds. From now on let $n\geq s$.

	\begin{itemize}
		\item[$\scriptstyle\blacklozenge$] \emph{Step~1: the inductive hypothesis.} By Lemma~\ref{lem:recursion},
		\begin{equation}\label{eq:upper-rec}
			D_{k+1}(n)=\sum_{q=0}^{\floor*{\frac{n-s}{b_{k+1}}}}D_k(n-s-qb_{k+1}),
		\end{equation}
		and every argument $n-s-qb_{k+1}$ with $0\leq q\leq\floor*{\frac{n-s}{b_{k+1}}}$ is a non-negative
		multiple of $d_k$. So we may apply the inductive hypothesis~\eqref{eq:upper-hyp} with
		$m=n-s-qb_{k+1}$: for $0\leq q\leq\floor*{\frac{n-s}{b_{k+1}}}$,
		\begin{equation*}
			D_k(n-s-qb_{k+1})\leq\frac{(n-s-qb_{k+1}+\nabla_k)^{k-1}}{(k-1)!\,b_2\dotsm b_k}.
		\end{equation*}
		Summing over $0\leq q\leq\floor*{\frac{n-s}{b_{k+1}}}$, by~\eqref{eq:upper-rec},
		\begin{equation}\label{eq:upper-step1}
			D_{k+1}(n)
			=\sum_{q=0}^{\floor*{\frac{n-s}{b_{k+1}}}}D_k(n-s-qb_{k+1})
			\leq\frac{1}{(k-1)!\,b_2\dotsm b_k}\sum_{q=0}^{\floor*{\frac{n-s}{b_{k+1}}}}(n-s-qb_{k+1}+\nabla_k)^{k-1}.
		\end{equation}

		\item[$\scriptstyle\blacklozenge$] \emph{Step~2: comparison with a progression.} The powers in~\eqref{eq:upper-step1} depend on
		$s$. We replace them by powers of the numbers $x-q$, where
		\begin{equation*}
			x\coloneqq\frac{n+\nabla_k}{b_{k+1}},
		\end{equation*}
		because the sums of such powers are estimated by Lemma~\ref{lem:sum-bounds}. Note that
		$x\geq0$, since $n\geq0$ and $\nabla_k\geq0$. Let $0\leq q\leq\floor*{\frac{n-s}{b_{k+1}}}$. Since
		$n-s-qb_{k+1}\geq0$ and $\nabla_k\geq0$ for the first inequality, since $s\geq0$ for the second
		one, and by the definition of $x$ for the equality,
		\begin{equation*}
			0
			\leq n-s-qb_{k+1}+\nabla_k
			\leq n-qb_{k+1}+\nabla_k
			=b_{k+1}\,(x-q).
		\end{equation*}
		Since $t\mapsto t^{k-1}$ is non-decreasing on $\R_{\geq0}$, it follows that
		\begin{equation*}
			(n-s-qb_{k+1}+\nabla_k)^{k-1}
			\leq\bigl(b_{k+1}\,(x-q)\bigr)^{k-1}
			=b_{k+1}^{k-1}\,(x-q)^{k-1}.
		\end{equation*}
		Inserting this estimate in~\eqref{eq:upper-step1}, we get
		\begin{equation}\label{eq:upper-step2}
			D_{k+1}(n)\leq\frac{b_{k+1}^{k-1}}{(k-1)!\,b_2\dotsm b_k}\sum_{q=0}^{\floor*{\frac{n-s}{b_{k+1}}}}(x-q)^{k-1}.
		\end{equation}

		\item[$\scriptstyle\blacklozenge$] \emph{Step~3: the sum of powers.} Lemma~\ref{lem:sum-bounds} concerns the sum over
		$0\leq q\leq\floor{x}$, while the sum in~\eqref{eq:upper-step2} stops at
		$\floor*{\frac{n-s}{b_{k+1}}}$. Since $s\geq0$ and $\nabla_k\geq0$,
		\begin{equation*}
			\frac{n-s}{b_{k+1}}
			\leq\frac{n}{b_{k+1}}
			\leq\frac{n+\nabla_k}{b_{k+1}}
			=x,
			\qquad\text{hence}\qquad
			\floor*{\frac{n-s}{b_{k+1}}}\leq\floor{x}.
		\end{equation*}
		Moreover, for every integer $q$ with $\floor*{\frac{n-s}{b_{k+1}}}<q\leq\floor{x}$ we have
		$q\leq x$, so $(x-q)^{k-1}\geq0$: extending the sum up to $\floor{x}$ only adds non-negative
		terms. Then Lemma~\ref{lem:sum-bounds} with $j=k-1\geq1$ applies:
		\begin{equation*}
			\sum_{q=0}^{\floor*{\frac{n-s}{b_{k+1}}}}(x-q)^{k-1}
			\leq\sum_{q=0}^{\floor{x}}(x-q)^{k-1}
			\leq\frac{\bigl(x+\frac12\bigr)^{k}}{k}.
		\end{equation*}
		Inserting this estimate in~\eqref{eq:upper-step2}, and using $(k-1)!\,k=k!$, we get
		\begin{equation}\label{eq:upper-step3}
			D_{k+1}(n)
			\leq\frac{b_{k+1}^{k-1}}{(k-1)!\,b_2\dotsm b_k}\cdot\frac{\bigl(x+\frac12\bigr)^{k}}{k}
			=\frac{b_{k+1}^{k-1}\bigl(x+\frac12\bigr)^{k}}{k!\,b_2\dotsm b_k}.
		\end{equation}

		\item[$\scriptstyle\blacklozenge$] \emph{Step~4: back to $n$.} By the definition of $x$, and then by
		Lemma~\ref{lem:basic}~(iv),
		\begin{equation*}
			b_{k+1}\Bigl(x+\frac12\Bigr)
			=n+\nabla_k+\frac{b_{k+1}}{2}
			=n+\nabla_{k+1},
			\qquad\text{that is}\qquad
			x+\frac12=\frac{n+\nabla_{k+1}}{b_{k+1}}.
		\end{equation*}
		Hence, by~\eqref{eq:upper-step3},
		\begin{equation*}
			D_{k+1}(n)
			\leq\frac{b_{k+1}^{k-1}\bigl(x+\frac12\bigr)^{k}}{k!\,b_2\dotsm b_k}
			=\frac{b_{k+1}^{k-1}\,(n+\nabla_{k+1})^{k}}{k!\,b_2\dotsm b_k\,b_{k+1}^{k}}
			=\frac{(n+\nabla_{k+1})^{k}}{k!\,b_2\dotsm b_k\,b_{k+1}}
			=\frac{(n+\nabla_{k+1})^{k}}{k!\,b_2\dotsm b_{k+1}},
		\end{equation*}
		which is~\eqref{eq:upper-goal}, that is, the statement for $k+1$.
	\end{itemize}
\end{proof}

\subsection{Lower bound}

\begin{thm}\label{thm:lower}
	Let $k\in\N_{\geq2}$ and let $n\in\N$ be a multiple of $d_k$ with $n\geq\Delta_k$. Then
	\begin{equation*}
		D_k(n)>\frac{(n-\Delta_k)^{k-1}}{(k-1)!\,b_2\dotsm b_k}.
	\end{equation*}
\end{thm}
\begin{proof}
	We prove the statement by induction on $k$.

	\emph{Base case $k=2$.} For $k=2$ the statement reads $D_2(n)>\frac{n-\Delta_2}{b_2}$ for every
	multiple $n\geq\Delta_2$ of $d_2$, which is the lower bound of Proposition~\ref{prop:two}.

	\emph{Inductive step.} Assume the statement for some $k\geq2$, that is, for every multiple
	$m\in\N$ of $d_k$ with $m\geq\Delta_k$,
	\begin{equation}\label{eq:lower-hyp}
		D_k(m)>\frac{(m-\Delta_k)^{k-1}}{(k-1)!\,b_2\dotsm b_k}.
	\end{equation}
	Let $n\in\N$ be a multiple of $d_{k+1}$ with $n\geq\Delta_{k+1}$, and let $s$ be the integer of
	Lemma~\ref{lem:recursion}. We have to prove that
	\begin{equation}\label{eq:lower-goal}
		D_{k+1}(n)>\frac{(n-\Delta_{k+1})^{k}}{k!\,b_2\dotsm b_{k+1}}.
	\end{equation}
	Throughout the inductive step let
	\begin{equation*}
		z\coloneqq\frac{n-\Delta_{k+1}}{b_{k+1}};
	\end{equation*}
	note that $z\geq0$, since $n\geq\Delta_{k+1}$.

	\begin{itemize}
		\item[$\scriptstyle\blacklozenge$] \emph{Step~1: the terms of the recursion that we keep.} The inductive
		hypothesis~\eqref{eq:lower-hyp} applies only to arguments $m\geq\Delta_k$. We show that the
		arguments $n-s-qb_{k+1}$ with $0\leq q\leq\floor{z}$ have this property, and that these terms
		belong to the sum of Lemma~\ref{lem:recursion}. This step does not use the inductive hypothesis.

		By Lemma~\ref{lem:recursion}, $s\leq b_{k+1}-a_{k+1}$; by Lemma~\ref{lem:basic}~(iv),
		$\Delta_k+b_{k+1}-a_{k+1}=\Delta_{k+1}$; and $n-\Delta_{k+1}=b_{k+1}z$ by the definition of $z$.
		Hence, for every integer $q$,
		\begin{equation}\label{eq:lower-arg}
			(n-s-qb_{k+1})-\Delta_k
			\geq n-(b_{k+1}-a_{k+1})-qb_{k+1}-\Delta_k
			=n-\Delta_{k+1}-qb_{k+1}
			=b_{k+1}\,(z-q).
		\end{equation}
		If $0\leq q\leq\floor{z}$, then $q\leq z$, so the right-hand side of~\eqref{eq:lower-arg} is
		non-negative: the argument $n-s-qb_{k+1}$ is at least $\Delta_k$.

		For $q=0$, \eqref{eq:lower-arg} reads $n-s-\Delta_k\geq b_{k+1}z$. Dividing by $b_{k+1}$, and
		using $\Delta_k\geq0$ (Lemma~\ref{lem:basic}~(iv)),
		\begin{equation}\label{eq:lower-floor}
			z
			\leq\frac{n-s-\Delta_k}{b_{k+1}}
			\leq\frac{n-s}{b_{k+1}},
			\qquad\text{hence}\qquad
			\floor{z}\leq\floor*{\frac{n-s}{b_{k+1}}} .
		\end{equation}
		So the terms with $0\leq q\leq\floor{z}$ are terms of the sum~\eqref{eq:recursion} of
		Lemma~\ref{lem:recursion}. The other terms of that sum are non-negative, because a denumerant is
		a number of solutions; dropping them can only decrease the sum. Therefore
		\begin{equation}\label{eq:lower-keep}
			D_{k+1}(n)
			=\sum_{q=0}^{\floor*{\frac{n-s}{b_{k+1}}}}D_k(n-s-qb_{k+1})
			\geq\sum_{q=0}^{\floor{z}}D_k(n-s-qb_{k+1}),
		\end{equation}
		where every argument on the right-hand side is a multiple of $d_k$
		(Lemma~\ref{lem:recursion}) and is at least $\Delta_k$.

		\item[$\scriptstyle\blacklozenge$] \emph{Step~2: the inductive hypothesis.} Let $0\leq q\leq\floor{z}$. By Step~1 we may apply the
		inductive hypothesis~\eqref{eq:lower-hyp} with $m=n-s-qb_{k+1}$:
		\begin{equation*}
			D_k(n-s-qb_{k+1})>\frac{(n-s-qb_{k+1}-\Delta_k)^{k-1}}{(k-1)!\,b_2\dotsm b_k}.
		\end{equation*}
		By~\eqref{eq:lower-arg}, $n-s-qb_{k+1}-\Delta_k\geq b_{k+1}\,(z-q)\geq0$, and $t\mapsto t^{k-1}$ is
		non-decreasing on $\R_{\geq0}$. Hence
		\begin{equation*}
			(n-s-qb_{k+1}-\Delta_k)^{k-1}
			\geq\bigl(b_{k+1}\,(z-q)\bigr)^{k-1}
			=b_{k+1}^{k-1}\,(z-q)^{k-1},
			\qquad\text{so that}\qquad
			D_k(n-s-qb_{k+1})>\frac{b_{k+1}^{k-1}\,(z-q)^{k-1}}{(k-1)!\,b_2\dotsm b_k}.
		\end{equation*}
		Summing over $0\leq q\leq\floor{z}$ (the sum contains at least the term $q=0$, so the
		inequality remains strict), by~\eqref{eq:lower-keep},
		\begin{equation}\label{eq:lower-step2}
			D_{k+1}(n)
			\geq\sum_{q=0}^{\floor{z}}D_k(n-s-qb_{k+1})
			>\frac{b_{k+1}^{k-1}}{(k-1)!\,b_2\dotsm b_k}\sum_{q=0}^{\floor{z}}(z-q)^{k-1}.
		\end{equation}

		\item[$\scriptstyle\blacklozenge$] \emph{Step~3: the sum of powers.} By the lower bound of Lemma~\ref{lem:sum-bounds} with
		$j=k-1\geq1$ and $x=z$, and since $\frac{z^{k-1}}{2}\geq0$,
		\begin{equation*}
			\sum_{q=0}^{\floor{z}}(z-q)^{k-1}
			\geq\frac{z^{k}}{k}+\frac{z^{k-1}}{2}
			\geq\frac{z^{k}}{k}.
		\end{equation*}
		Inserting this estimate in~\eqref{eq:lower-step2}, and using $(k-1)!\,k=k!$, we get
		\begin{equation}\label{eq:lower-step3}
			D_{k+1}(n)
			>\frac{b_{k+1}^{k-1}}{(k-1)!\,b_2\dotsm b_k}\cdot\frac{z^{k}}{k}
			=\frac{b_{k+1}^{k-1}z^{k}}{k!\,b_2\dotsm b_k}.
		\end{equation}

		\item[$\scriptstyle\blacklozenge$] \emph{Step~4: back to $n$.} By the definition of $z$, $b_{k+1}z=n-\Delta_{k+1}$, that is
		$z=\frac{n-\Delta_{k+1}}{b_{k+1}}$. Hence, by~\eqref{eq:lower-step3},
		\begin{equation*}
			D_{k+1}(n)
			>\frac{b_{k+1}^{k-1}z^{k}}{k!\,b_2\dotsm b_k}
			=\frac{b_{k+1}^{k-1}\,(n-\Delta_{k+1})^{k}}{k!\,b_2\dotsm b_k\,b_{k+1}^{k}}
			=\frac{(n-\Delta_{k+1})^{k}}{k!\,b_2\dotsm b_k\,b_{k+1}}
			=\frac{(n-\Delta_{k+1})^{k}}{k!\,b_2\dotsm b_{k+1}},
		\end{equation*}
		which is~\eqref{eq:lower-goal}, that is, the statement for $k+1$.
	\end{itemize}
\end{proof}

\begin{cor}\label{cor:positive}
	Let $k\in\N_{\geq2}$ and let $n\in\N$ be a multiple of $d_k$ with $n\geq\Delta_k$. Then
	$D_k(n)\geq1$. In particular, if $d_k=1$, every integer $n\geq\Delta_k$ is representable, and
	$g(a_1,\dotsc,a_k)<\Delta_k$.
\end{cor}
\begin{proof}
	By Theorem~\ref{thm:lower}, $D_k(n)$ is an integer strictly greater than a non-negative number.
\end{proof}

\begin{proof}[Proof of Theorem~\ref{thm:main}]
	By Lemma~\ref{lem:basic}~(iii), $\frac{1}{b_2\dotsm b_k}=\frac{d_k}{a_1\dotsm a_k}$; the upper
	bound is Theorem~\ref{thm:upper} and the lower bound is Theorem~\ref{thm:lower}.
\end{proof}

\subsection{A sharper lower bound}

\begin{thm}\label{thm:lower-poly}
	Let $k\in\N_{\geq2}$ and let $n\in\N$ be a multiple of $d_k$ with $n\geq\Delta_k$. Then
	\begin{equation*}
		D_k(n)>\frac{1}{b_2\dotsm b_k}\sum_{i=0}^{k-2}\bbinom{k-2}{i}\frac{(n-\Delta_k)^{k-1-i}}{(k-1-i)!}.
	\end{equation*}
\end{thm}
\begin{proof}
	For $k\in\N_{\geq2}$ and $t\in\R$ let
	\begin{equation}\label{eq:P-def}
		P_k(t)\coloneqq\sum_{i=0}^{k-2}\bbinom{k-2}{i}\frac{t^{k-1-i}}{(k-1-i)!}.
	\end{equation}
	The sum in the statement is $P_k(n-\Delta_k)$. So the statement reads: for every $k\geq2$ and
	every multiple $n\in\N$ of $d_k$ with $n\geq\Delta_k$,
	\begin{equation}\label{eq:P-claim}
		D_k(n)>\frac{P_k(n-\Delta_k)}{b_2\dotsm b_k}.
	\end{equation}
	We shall use that $P_k$ is non-decreasing on $\R_{\geq0}$: in~\eqref{eq:P-def} the
	coefficients $\bbinom{k-2}{i}$ are non-negative (Lemma~\ref{lem:sym-poly}~(i)), and every
	exponent $k-1-i$ is at least $1$, so every term is non-decreasing on $\R_{\geq0}$.

	We prove~\eqref{eq:P-claim} by induction on $k$.

	\emph{Base case $k=2$.} For $k=2$ the sum in~\eqref{eq:P-def} has the single term $i=0$, and
	$\bbinom{0}{0}=1$ (Lemma~\ref{lem:sym-poly}~(i)), so $P_2(t)=t$. Hence~\eqref{eq:P-claim} reads
	$D_2(n)>\frac{n-\Delta_2}{b_2}$ for every multiple $n\geq\Delta_2$ of $d_2$, which is the lower
	bound of Proposition~\ref{prop:two}.

	\emph{Inductive step.} Assume~\eqref{eq:P-claim} for some $k\geq2$, that is, for every multiple
	$m\in\N$ of $d_k$ with $m\geq\Delta_k$,
	\begin{equation}\label{eq:P-hyp}
		D_k(m)>\frac{P_k(m-\Delta_k)}{b_2\dotsm b_k}.
	\end{equation}
	Let $n\in\N$ be a multiple of $d_{k+1}$ with $n\geq\Delta_{k+1}$, and let $s$ be the integer of
	Lemma~\ref{lem:recursion}. We have to prove that
	\begin{equation}\label{eq:P-goal}
		D_{k+1}(n)>\frac{P_{k+1}(n-\Delta_{k+1})}{b_2\dotsm b_{k+1}}.
	\end{equation}
	As in the proof of Theorem~\ref{thm:lower}, let
	\begin{equation*}
		z\coloneqq\frac{n-\Delta_{k+1}}{b_{k+1}}\geq0 .
	\end{equation*}

	\begin{itemize}
		\item[$\scriptstyle\blacklozenge$] \emph{Step~1: the terms of the recursion that we keep.} Step~1 of the proof of
		Theorem~\ref{thm:lower} concerns the same $n$, $s$ and $z$, and it does not use the inductive
		hypothesis of that theorem: it uses only Lemmas~\ref{lem:basic} and~\ref{lem:recursion}. So its
		two conclusions hold here. They are~\eqref{eq:lower-arg}, which for $0\leq q\leq\floor{z}$ gives
		\begin{equation}\label{eq:P-arg}
			n-s-qb_{k+1}-\Delta_k\geq b_{k+1}\,(z-q)\geq0,
		\end{equation}
		and~\eqref{eq:lower-keep}, that is
		\begin{equation}\label{eq:P-keep}
			D_{k+1}(n)\geq\sum_{q=0}^{\floor{z}}D_k(n-s-qb_{k+1}),
		\end{equation}
		where every argument is a multiple of $d_k$ and is at least $\Delta_k$.

		\item[$\scriptstyle\blacklozenge$] \emph{Step~2: the inductive hypothesis.} Let $0\leq q\leq\floor{z}$. By Step~1 we may apply the
		inductive hypothesis~\eqref{eq:P-hyp} with $m=n-s-qb_{k+1}$; then, by~\eqref{eq:P-arg} and
		because $P_k$ is non-decreasing on $\R_{\geq0}$,
		\begin{equation*}
			D_k(n-s-qb_{k+1})
			>\frac{P_k(n-s-qb_{k+1}-\Delta_k)}{b_2\dotsm b_k}
			\geq\frac{P_k\bigl(b_{k+1}\,(z-q)\bigr)}{b_2\dotsm b_k}.
		\end{equation*}
		Summing over $0\leq q\leq\floor{z}$ (the sum contains at least the term $q=0$, so the
		inequality remains strict), by~\eqref{eq:P-keep},
		\begin{equation}\label{eq:P-step2}
			D_{k+1}(n)
			\geq\sum_{q=0}^{\floor{z}}D_k(n-s-qb_{k+1})
			>\frac{1}{b_2\dotsm b_k}\sum_{q=0}^{\floor{z}}P_k\bigl(b_{k+1}\,(z-q)\bigr).
		\end{equation}

		\item[$\scriptstyle\blacklozenge$] \emph{Step~3: the sums of powers.} By the definition~\eqref{eq:P-def} of $P_k$, and exchanging
		the two finite sums,
		\begin{equation}\label{eq:P-exchange}
			\sum_{q=0}^{\floor{z}}P_k\bigl(b_{k+1}\,(z-q)\bigr)
			=\sum_{q=0}^{\floor{z}}\sum_{i=0}^{k-2}\bbinom{k-2}{i}\frac{b_{k+1}^{k-1-i}\,(z-q)^{k-1-i}}{(k-1-i)!}
			=\sum_{i=0}^{k-2}\bbinom{k-2}{i}\frac{b_{k+1}^{k-1-i}}{(k-1-i)!}\sum_{q=0}^{\floor{z}}(z-q)^{k-1-i}.
		\end{equation}
		Let $0\leq i\leq k-2$. The exponent $k-1-i$ is at least $1$, so the lower bound of
		Lemma~\ref{lem:sum-bounds} with $j=k-1-i$ and $x=z$ applies to the inner sum:
		\begin{equation*}
			\sum_{q=0}^{\floor{z}}(z-q)^{k-1-i}\geq\frac{z^{k-i}}{k-i}+\frac{z^{k-1-i}}{2}.
		\end{equation*}
		Multiplying by $\frac{b_{k+1}^{k-1-i}}{(k-1-i)!}>0$, and using $(k-1-i)!\,(k-i)=(k-i)!$,
		\begin{equation*}
			\frac{b_{k+1}^{k-1-i}}{(k-1-i)!}\sum_{q=0}^{\floor{z}}(z-q)^{k-1-i}
			\geq\frac{b_{k+1}^{k-1-i}z^{k-i}}{(k-i)!}+\frac{b_{k+1}^{k-1-i}z^{k-1-i}}{2\,(k-1-i)!}
			=\frac{1}{b_{k+1}}\biggl(\frac{(b_{k+1}z)^{k-i}}{(k-i)!}+\frac{b_{k+1}}{2}\cdot\frac{(b_{k+1}z)^{k-1-i}}{(k-1-i)!}\biggr).
		\end{equation*}
		The coefficients $\bbinom{k-2}{i}$ are non-negative (Lemma~\ref{lem:sym-poly}~(i)), so we may
		insert this estimate in~\eqref{eq:P-exchange}:
		\begin{equation}\label{eq:P-step3}
			\sum_{q=0}^{\floor{z}}P_k\bigl(b_{k+1}\,(z-q)\bigr)
			\geq\frac{1}{b_{k+1}}\Biggl(\,\sum_{i=0}^{k-2}\bbinom{k-2}{i}\frac{(b_{k+1}z)^{k-i}}{(k-i)!}
			+\frac{b_{k+1}}{2}\sum_{i=0}^{k-2}\bbinom{k-2}{i}\frac{(b_{k+1}z)^{k-1-i}}{(k-1-i)!}\Biggr).
		\end{equation}

		\item[$\scriptstyle\blacklozenge$] \emph{Step~4: collecting the terms.} We show that the expression in brackets
		in~\eqref{eq:P-step3} is $P_{k+1}(b_{k+1}z)$. In its second sum we replace $i$ by $i-1$:
		\begin{equation*}
			\sum_{i=0}^{k-2}\bbinom{k-2}{i}\frac{(b_{k+1}z)^{k-1-i}}{(k-1-i)!}
			=\sum_{i=1}^{k-1}\bbinom{k-2}{i-1}\frac{(b_{k+1}z)^{k-i}}{(k-i)!}.
		\end{equation*}
		By Lemma~\ref{lem:sym-poly}~(i), $\bbinom{k-2}{-1}=0$ and $\bbinom{k-2}{k-1}=0$. So we may add
		the term $i=0$ to this sum and the term $i=k-1$ to the first sum in~\eqref{eq:P-step3}, since
		both are $0$; then both sums run over $0\leq i\leq k-1$, and the expression in brackets is
		\begin{equation*}
			\sum_{i=0}^{k-1}\Biggl(\bbinom{k-2}{i}+\frac{b_{k+1}}{2}\,\bbinom{k-2}{i-1}\Biggr)\frac{(b_{k+1}z)^{k-i}}{(k-i)!}
			=\sum_{i=0}^{k-1}\bbinom{k-1}{i}\frac{(b_{k+1}z)^{k-i}}{(k-i)!}
			=P_{k+1}(b_{k+1}z),
		\end{equation*}
		where the first equality is Lemma~\ref{lem:sym-poly}~(ii) with $m=k-1$, and the second one is
		the definition~\eqref{eq:P-def} with $k+1$ in place of $k$. Therefore~\eqref{eq:P-step3} reads
		\begin{equation*}
			\sum_{q=0}^{\floor{z}}P_k\bigl(b_{k+1}\,(z-q)\bigr)\geq\frac{P_{k+1}(b_{k+1}z)}{b_{k+1}}.
		\end{equation*}
		Inserting this in~\eqref{eq:P-step2}, and using $b_{k+1}z=n-\Delta_{k+1}$ and
		$b_2\dotsm b_k\,b_{k+1}=b_2\dotsm b_{k+1}$,
		\begin{equation*}
			D_{k+1}(n)
			>\frac{1}{b_2\dotsm b_k}\cdot\frac{P_{k+1}(b_{k+1}z)}{b_{k+1}}
			=\frac{P_{k+1}(n-\Delta_{k+1})}{b_2\dotsm b_{k+1}},
		\end{equation*}
		which is~\eqref{eq:P-goal}, that is, \eqref{eq:P-claim} for $k+1$.
	\end{itemize}
\end{proof}

\begin{rem}\label{rem:lower-compare}
	The term $i=0$ of Theorem~\ref{thm:lower-poly} is the bound of Theorem~\ref{thm:lower},
	since $\bbinom{k-2}{0}=1$, and the other terms are non-negative. Hence
	Theorem~\ref{thm:lower-poly} is at least as sharp as Theorem~\ref{thm:lower}, and strictly
	sharper when $k\geq3$ and $n>\Delta_k$, since then the term $i=1$ is positive.
\end{rem}

% =====================================================================
\section{Asymptotic behaviour}\label{s:asymptotic}
% =====================================================================
The two bounds of Theorem~\ref{thm:main} have the same leading term, and each of them differs
from it by $O\bigl(n^{k-2}\bigr)$; so they give at once the asymptotic behaviour of $D_k(n)$
stated in the next corollary, which for $d_k=1$ is Nathanson's formula~\eqref{eq:intro-nathanson}.
We omit this routine verification, because it adds nothing. The purpose of this section is a
different one: to show that the technique of Section~\ref{s:main}, in a simplified form, solves
this similar but weaker problem directly, and can be used in place of the argument of
Nathanson~\cite{nathanson}.

\begin{cor}\label{cor:asymptotic}
	Let $k\in\N_{\geq2}$, and set $c_k\coloneqq\frac{d_k}{(k-1)!\,a_1\dotsm a_k}$. Then, along the
	multiples of $d_k$,
	\begin{equation*}
		D_k(n)=c_k\,n^{k-1}+O\bigl(n^{k-2}\bigr),
	\end{equation*}
	that is, there are $C\in\R_{>0}$ and $n_0\in\N$ such that
	$\abs{D_k(n)-c_k\,n^{k-1}}\leq C\,n^{k-2}$ for every multiple $n\geq n_0$ of $d_k$. In
	particular, if $d_k=1$, the formula holds for $n\in\N$.
\end{cor}

\begin{rem}\label{rem:asymptotic-dk}
	If $d_k>1$, then $D_k(n)=0$ for every $n$ that is not a multiple of $d_k$ (Lemma~\ref{lem:basic}~(vi)),
	so no asymptotic formula can hold for all $n$; Corollary~\ref{cor:asymptotic} describes the
	behaviour along the multiples of $d_k$, where the leading term carries the factor $d_k$.
\end{rem}

We do not write the simplified proof in full, because its computations are those of the proofs
of Theorems~\ref{thm:upper} and~\ref{thm:lower}, with every explicit constant replaced by an
$O$-term; we describe it, and then we compare it with the proof in~\cite{nathanson}.

\emph{The simplified argument.} One proves by induction on $k\in\N_{\geq2}$ that, along the
multiples of $d_k$,
\begin{equation}\label{eq:asymptotic-ind}
	D_k(n)=\frac{n^{k-1}}{(k-1)!\,b_2\dotsm b_k}+O\bigl(n^{k-2}\bigr),
\end{equation}
which by Lemma~\ref{lem:basic}~(iii) is the formula of Corollary~\ref{cor:asymptotic}. The base
case $k=2$ is the formula for $D_2(n)$ of Proposition~\ref{prop:two}, because
$\floor*{\frac{n-s}{b_2}}+1=\frac{n}{b_2}+O(1)$. In the inductive step, let $n$ be a multiple of
$d_{k+1}$, let $s$ be as in Lemma~\ref{lem:recursion}, and set $y\coloneqq\frac{n-s}{b_{k+1}}$,
so that $n-s-qb_{k+1}=b_{k+1}\,(y-q)$. Then
\begin{equation}\label{eq:asymptotic-step}
	\begin{aligned}
		D_{k+1}(n)
		&=\sum_{q=0}^{\floor{y}}D_k\bigl(n-s-qb_{k+1}\bigr)
		=\frac{b_{k+1}^{k-1}}{(k-1)!\,b_2\dotsm b_k}\sum_{q=0}^{\floor{y}}(y-q)^{k-1}+O\bigl(n^{k-1}\bigr)\\
		&=\frac{b_{k+1}^{k-1}\,y^{k}}{k!\,b_2\dotsm b_k}+O\bigl(n^{k-1}\bigr)
		=\frac{n^{k}}{k!\,b_2\dotsm b_{k+1}}+O\bigl(n^{k-1}\bigr),
	\end{aligned}
\end{equation}
which is~\eqref{eq:asymptotic-ind} for $k+1$. The first equality is the
recursion~\eqref{eq:recursion}. The second one is the inductive hypothesis, used in each of the
at most $\frac{n}{b_{k+1}}+1$ terms, each with an error $O\bigl(n^{k-2}\bigr)$. The third one is
Lemma~\ref{lem:sum-bounds}, of which only the consequence
$\sum_{q=0}^{\floor{y}}(y-q)^{j}=\frac{y^{j+1}}{j+1}+O\bigl(y^{j}\bigr)$ is needed. The fourth one
is $b_{k+1}^{k-1}\,y^{k}=\frac{(n-s)^{k}}{b_{k+1}}=\frac{n^{k}}{b_{k+1}}+O\bigl(n^{k-1}\bigr)$,
because $0\leq s<b_{k+1}$. These are the steps of the proofs of Theorems~\ref{thm:upper}
and~\ref{thm:lower}, stripped of everything whose only role there is to keep the constants
explicit: the shifts
$\nabla_k$ and $\Delta_k$, the terms added to the sum or dropped from it, and the two-sided
estimate of Lemma~\ref{lem:sum-bounds}.

\emph{Comparison with the proof of Nathanson.} The proof in~\cite{nathanson} is also an
induction on $k$, and it rests on the same arithmetic fact: for $d_{k+1}=1$ the admissible
values of $x_{k+1}$ are shown to be $u+\ell d_k$, with a fixed integer $0\leq u<d_k$ and
$\ell\in\N$; since then $b_{k+1}=d_ka_{k+1}$ (Lemma~\ref{lem:basic}~(ii)), we have
$a_{k+1}x_{k+1}=a_{k+1}u+\ell\,b_{k+1}$, the progression $s+qb_{k+1}$ of
Lemma~\ref{lem:recursion}. The two proofs differ in two respects.
\begin{enumerate}
	\item[\textup{(i)}] In~\cite{nathanson} the inductive hypothesis is stated for coprime
	      integers only, so $a_1,\dotsc,a_k$ and the argument of $D_k$ are first divided by
	      $d_k$. Here~\eqref{eq:asymptotic-ind} holds for arbitrary $a_1,\dotsc,a_k$ along the
	      multiples of $d_k$, and applies directly to the terms of the recursion.
	\item[\textup{(ii)}] In~\cite{nathanson} every power in the sum is expanded by the binomial
	      theorem, the sums $\sum_{\ell=0}^{r}\ell^j$ are replaced by
	      $\frac{r^{j+1}}{j+1}+O\bigl(r^j\bigr)$, and the leading coefficient is computed through
	      an identity for alternating sums of binomial coefficients. Here the sum
	      in~\eqref{eq:asymptotic-step} is not expanded: it is estimated as a whole.
\end{enumerate}

% =====================================================================
\section{Comparison with the known results}\label{s:comparison}
% =====================================================================
We compare the bounds of Section~\ref{s:main} with the results recalled in the Introduction.

\subsection{The case \texorpdfstring{$a_1=1$}{a1=1}}
Let $a_1=1$. By Lemma~\ref{lem:basic}~(v), $d_k=1$, $b_i=a_i$ for every $i$, and $\Delta_k=0$.
Consequently:
\begin{itemize}
	\item every $n\in\N$ is a multiple of $d_k$ and satisfies $n\geq\Delta_k$;
	\item $b_2\dotsm b_k=a_2\dotsm a_k=a_1a_2\dotsm a_k$;
	\item $\nabla_k=b_2+\frac12\sum_{i=3}^{k}b_i=a_2+\frac12\sum_{i=3}^{k}a_i=s_k$, the constant
	      of~\eqref{eq:intro-blom-froberg};
	\item $\bbinom{k-2}{i}=e_i\bigl(\frac{b_3}{2},\dotsc,\frac{b_k}{2}\bigr)
	      =e_i\bigl(\frac{a_3}{2},\dotsc,\frac{a_k}{2}\bigr)$.
\end{itemize}
Hence Theorems~\ref{thm:upper} and~\ref{thm:lower-poly} give, for every $n\in\N$,
\begin{equation*}
	\frac{1}{a_1\dotsm a_k}\sum_{i=0}^{k-2}\bbinom{k-2}{i}\frac{n^{k-1-i}}{(k-1-i)!}
	<D_k(n)\leq\frac{(n+s_k)^{k-1}}{(k-1)!\,a_1\dotsm a_k}.
\end{equation*}
These are exactly the two inequalities of Blom and Fr\"oberg~\cite{blom-froberg}. The upper
bound, together with the term $i=0$ of the lower bound, is their inequality~(A) in Section~2.3
of~\cite{blom-froberg}, that is~\eqref{eq:intro-blom-froberg}. The whole lower bound is their
inequality~(B) in Section~2.4 of~\cite{blom-froberg}: their coefficients $S_{i,k}$, the sums of the
products of $i$ distinct numbers among $\frac{a_3}{2},\dotsc,\frac{a_k}{2}$, are precisely
$e_i\bigl(\frac{a_3}{2},\dotsc,\frac{a_k}{2}\bigr)=\bbinom{k-2}{i}$. Both their lower bounds are
stated with $\leq$; here they hold with strict inequality.

\subsection{The bounds of Agnarsson}
We first write the constants $A_k$ and $B_k$ of~\eqref{eq:intro-AB} in our notation. We have
$\frac{d_{k-1}}{d_k}\,a_k=b_k$ by Lemma~\ref{lem:basic}~(ii), so the increments
in~\eqref{eq:intro-AB} are $b_k$ and $b_k-a_k-1$. Since $A_1=B_1=0$, summing them we get
\begin{equation}\label{eq:AB-ours}
	A_k=\sum_{i=2}^{k}b_i,
	\qquad
	B_k
	=\sum_{i=2}^{k}(b_i-a_i-1)
	=\sum_{i=2}^{k}(b_i-a_i)-(k-1)
	=\Delta_k-(k-1).
\end{equation}
By the first of these formulas and the definition of $\nabla_k$,
\begin{equation}\label{eq:A-nabla}
	A_k-\nabla_k
	=\biggl(b_2+\sum_{i=3}^{k}b_i\biggr)-\biggl(b_2+\frac12\sum_{i=3}^{k}b_i\biggr)
	=\frac12\sum_{i=3}^{k}b_i .
\end{equation}
Finally, Agnarsson's bounds~\eqref{eq:intro-agnarsson} and ours have the same factor, because
$\frac{d_k}{a_1\dotsm a_k}=\frac{1}{b_2\dotsm b_k}$ by Lemma~\ref{lem:basic}~(iii); so in the
comparisons below we leave it out.

\begin{prop}[Agnarsson's lower bound and Theorem~\ref{thm:lower}]\label{prop:compare-lower}
	Let $k\in\N_{\geq2}$. For every $n\geq\Delta_k$, the lower bound of Agnarsson is greater than
	the lower bound of Theorem~\ref{thm:lower}.
\end{prop}
\begin{proof}
	Let $n\geq\Delta_k$. We have to prove that
	\begin{equation*}
		\binom{n-B_k}{k-1}>\frac{(n-\Delta_k)^{k-1}}{(k-1)!}.
	\end{equation*}
	By~\eqref{eq:AB-ours}, $n-B_k=n-\Delta_k+k-1$. So the numerator of the binomial coefficient
	is the product of the $k-1$ consecutive integers from $n-\Delta_k+k-1$ down to
	$n-\Delta_k+1$. Each of them is at least $n-\Delta_k+1$, which is positive, so
	\begin{equation*}
		\binom{n-B_k}{k-1}
		=\binom{n-\Delta_k+k-1}{k-1}
		=\frac{(n-\Delta_k+k-1)\dotsm(n-\Delta_k+1)}{(k-1)!}
		\geq\frac{(n-\Delta_k+1)^{k-1}}{(k-1)!}
		>\frac{(n-\Delta_k)^{k-1}}{(k-1)!},
	\end{equation*}
	where the last inequality holds because $0\leq n-\Delta_k<n-\Delta_k+1$ and
	$t\mapsto t^{k-1}$ is strictly increasing on $\R_{\geq0}$ for $k\geq2$.
\end{proof}

\begin{prop}[Agnarsson's upper bound and Theorem~\ref{thm:upper}]\label{prop:compare-upper}
	Let $k\in\N_{\geq2}$ with $\frac12\sum_{i=3}^{k}b_i\geq k-2$. For every $n\in\N$, the upper
	bound of Theorem~\ref{thm:upper} is less than or equal to the upper bound of Agnarsson.
\end{prop}
\begin{proof}
	Let $n\in\N$. We have to prove that
	\begin{equation*}
		\frac{(n+\nabla_k)^{k-1}}{(k-1)!}\leq\binom{n+A_k}{k-1}.
	\end{equation*}
	By~\eqref{eq:A-nabla}, the assumption is $A_k-\nabla_k\geq k-2$, that is,
	$n+A_k-k+2\geq n+\nabla_k$. The $k-1$
	factors of the numerator of the binomial coefficient are
	$n+A_k,\,n+A_k-1,\,\dotsc,\,n+A_k-k+2$. Each of them is at least $n+A_k-k+2$, which is at
	least $n+\nabla_k>0$, so
	\begin{equation*}
		\binom{n+A_k}{k-1}
		=\frac{(n+A_k)(n+A_k-1)\dotsm(n+A_k-k+2)}{(k-1)!}
		\geq\frac{(n+A_k-k+2)^{k-1}}{(k-1)!}
		\geq\frac{(n+\nabla_k)^{k-1}}{(k-1)!},
	\end{equation*}
	where the last inequality holds because $t\mapsto t^{k-1}$ is non-decreasing on
	$\R_{\geq0}$.
\end{proof}

\begin{prop}[Agnarsson's lower bound and Theorem~\ref{thm:lower-poly}]\label{prop:compare-lower-poly}
	Let $k\in\N_{\geq3}$ with $b_3+\dotsb+b_k>k$. For every sufficiently large $n$, the lower
	bound of Theorem~\ref{thm:lower-poly} is greater than the lower bound of Agnarsson.
\end{prop}
\begin{proof}
	We have to prove that, for every sufficiently large $n$,
	\begin{equation*}
		\sum_{i=0}^{k-2}\bbinom{k-2}{i}\frac{(n-\Delta_k)^{k-1-i}}{(k-1-i)!}>\binom{n-B_k}{k-1}.
	\end{equation*}
	Set $N\coloneqq n-\Delta_k$; by~\eqref{eq:AB-ours}, $n-B_k=N+k-1$. We compare the two
	leading terms of the two sides. On the left, since $\bbinom{k-2}{0}=1$ and
	$\bbinom{k-2}{1}=\frac{b_3+\dotsb+b_k}{2}$,
	\begin{equation*}
		\sum_{i=0}^{k-2}\bbinom{k-2}{i}\frac{N^{k-1-i}}{(k-1-i)!}
		=\frac{N^{k-1}}{(k-1)!}+\frac{b_3+\dotsb+b_k}{2\,(k-2)!}\,N^{k-2}+O\bigl(N^{k-3}\bigr).
	\end{equation*}
	On the right, expanding the product,
	\begin{align*}
		\binom{n-B_k}{k-1}
		&=\binom{N+k-1}{k-1}
		=\frac{(N+1)(N+2)\dotsm(N+k-1)}{(k-1)!}\\
		&=\frac{N^{k-1}}{(k-1)!}+\frac{1+2+\dotsb+(k-1)}{(k-1)!}\,N^{k-2}+O\bigl(N^{k-3}\bigr)
		=\frac{N^{k-1}}{(k-1)!}+\frac{k}{2\,(k-2)!}\,N^{k-2}+O\bigl(N^{k-3}\bigr).
	\end{align*}
	Hence the difference between the left and the right side is
	\begin{equation*}
		\frac{b_3+\dotsb+b_k-k}{2\,(k-2)!}\,N^{k-2}+O\bigl(N^{k-3}\bigr).
	\end{equation*}
	Since $b_3+\dotsb+b_k>k$, the coefficient of $N^{k-2}$ is positive, so the difference is
	positive for large $N$.
\end{proof}

\begin{rem}\label{rem:compare-agnarsson}
	\begin{enumerate}
		\item[\textup{(i)}] By~\eqref{eq:AB-ours}, the range $n\geq B_k+k-1$ of Agnarsson's
		      lower bound is $n\geq\Delta_k$, the range of Theorem~\ref{thm:lower}; so
		      Proposition~\ref{prop:compare-lower} covers the whole range of both bounds. Its
		      proof shows more: Agnarsson's lower bound is at least the bound of
		      Theorem~\ref{thm:lower} with $\Delta_k$ replaced by $\Delta_k-1$.
		\item[\textup{(ii)}] The assumption of Proposition~\ref{prop:compare-upper} holds if
		      $b_3,\dotsc,b_k\geq2$, since then
		      $\frac12\sum_{i=3}^{k}b_i\geq\frac12\cdot2\,(k-2)=k-2$. For $k=2$ it holds
		      trivially, and the two upper bounds coincide: $A_2=b_2=\nabla_2$ and
		      $\binom{n+b_2}{1}=n+b_2$.
		\item[\textup{(iii)}] In Proposition~\ref{prop:compare-lower-poly} the comparison may go
		      the other way for small $n$: at $n=\Delta_k$ every term of the sum in
		      Theorem~\ref{thm:lower-poly} contains a positive power of $n-\Delta_k=0$, so the
		      bound vanishes, while Agnarsson's equals $\frac{1}{b_2\dotsm b_k}$.
	\end{enumerate}
\end{rem}

\subsection{Positivity and the Frobenius number}
Our proof of the lower bound does not assume $D_k(n)>0$: positivity for $n\geq\Delta_k$ is
obtained along the induction (Corollary~\ref{cor:positive}). In particular, the assumption
$D_k(n)>0$ of Theorem~3.4 in~\cite{agnarsson} is always satisfied by the multiples
$n\geq\Delta_k$ of $d_k$, so that Agnarsson's lower bound holds on that range without further
assumptions.

For $d_k=1$, Corollary~\ref{cor:positive} says that every $n\geq\Delta_k$ is representable,
that is, $g(a_1,\dotsc,a_k)<\Delta_k$. Let us compare this with Brauer's
bound~\eqref{eq:intro-brauer}. By Lemma~\ref{lem:basic}~(ii), and separating the term $i=1$ of
the second sum, its right-hand side is
\begin{equation*}
	\sum_{i=2}^{k}\frac{d_{i-1}}{d_i}\,a_i-\sum_{i=1}^{k}a_i
	=\sum_{i=2}^{k}b_i-\sum_{i=1}^{k}a_i
	=\sum_{i=2}^{k}(b_i-a_i)-a_1
	=\Delta_k-a_1,
\end{equation*}
so Brauer's bound is $g(a_1,\dotsc,a_k)\leq\Delta_k-a_1$. Since $a_1\geq1$, as a bound on the
Frobenius number ours is weaker; the point is that positivity and the quantitative estimate
come out of the same elementary argument.

\subsection{Asymptotic estimates}
Gajdzica's estimate~\eqref{eq:intro-gajdzica} has the explicit constant $1$ in the error term,
while Corollary~\ref{cor:asymptotic} has an unspecified one. On the other hand,
\eqref{eq:intro-gajdzica} assumes the $a_i$ distinct and the two smallest of them coprime,
whereas Corollary~\ref{cor:asymptotic} holds for arbitrary $a_1,\dotsc,a_k$ along the multiples
of $d_k$.

\subsection{The order of the \texorpdfstring{$a_i$}{ai}}
$D_k(n)$ does not depend on the order of $a_1,\dotsc,a_k$, while $\Delta_k$ and $\nabla_k$ do
(as do Agnarsson's $A_k$ and $B_k$, see~\cite{agnarsson}): each bound may be applied to the most
favourable ordering.

% =====================================================================
\section{Final remarks}\label{s:final}
% =====================================================================

\subsection{Refined bounds at the last step}
In the proofs of Theorems~\ref{thm:upper} and~\ref{thm:lower}, the integer $s$ of
Lemma~\ref{lem:recursion} is estimated by $0\leq s\leq b_{k+1}-a_{k+1}$, uniformly in $n$.
Keeping $s$ at the last step of the induction gives sharper bounds, which depend on $n$ through $s$.

\begin{prop}\label{prop:refined}
	Let $k\in\N_{\geq2}$, let $n\in\N$ be a multiple of $d_k$, and let $s$ be the integer of
	Lemma~\ref{lem:recursion} for the passage from $k-1$ to $k$, that is, the integer with
	$0\leq s<b_k$, $a_k\mid s$ and $d_{k-1}\mid n-s$. For $k=2$ we read $\Delta_1$ as the empty sum
	$\sum_{i=2}^{1}(b_i-a_i)=0$.
	\begin{enumerate}
		\item[\textup{(i)}] If $n\geq s$, then
		      \begin{equation*}
		      	D_k(n)\leq\frac{(n-s+\nabla_k)^{k-1}}{(k-1)!\,b_2\dotsm b_k}.
		      \end{equation*}
		\item[\textup{(ii)}] If $n\geq s+\Delta_{k-1}$, then
		      \begin{equation*}
		      	D_k(n)>\frac{(n-s-\Delta_{k-1})^{k-1}}{(k-1)!\,b_2\dotsm b_k}.
		      \end{equation*}
	\end{enumerate}
\end{prop}
\begin{proof}
	\emph{The case $k=2$.} Here $\Delta_1=0$, so both (i) and (ii) assume $n\geq s$, and together
	they read
	\begin{equation*}
		\frac{n-s}{b_2}<D_2(n)\leq\frac{n-s+\nabla_2}{b_2}.
	\end{equation*}
	Since $n\geq s$, Proposition~\ref{prop:two} gives $D_2(n)=\floor*{\frac{n-s}{b_2}}+1$. For
	every $y\in\R$ we have $y<\floor{y}+1\leq y+1$; with $y=\frac{n-s}{b_2}$, and then by
	$\nabla_2=b_2$ (Lemma~\ref{lem:basic}~(iv)),
	\begin{equation*}
		\frac{n-s}{b_2}
		<D_2(n)
		\leq\frac{n-s}{b_2}+1
		=\frac{n-s+b_2}{b_2}
		=\frac{n-s+\nabla_2}{b_2}.
	\end{equation*}

	\emph{The case $k\geq3$.} By Lemma~\ref{lem:recursion}, for the passage from $k-1$ to $k$,
	\begin{equation}\label{eq:refined-rec}
		D_k(n)=\sum_{q=0}^{\floor*{\frac{n-s}{b_k}}}D_{k-1}(n-s-qb_k),
	\end{equation}
	where every argument $n-s-qb_k$ is a non-negative multiple of $d_{k-1}$. By
	Theorems~\ref{thm:lower} and~\ref{thm:upper} for $D_{k-1}$, we know that, for every multiple
	$m$ of $d_{k-1}$,
	\begin{equation}\label{eq:refined-known}
		\frac{(m-\Delta_{k-1})^{k-2}}{(k-2)!\,b_2\dotsm b_{k-1}}
		<D_{k-1}(m)
		\leq\frac{(m+\nabla_{k-1})^{k-2}}{(k-2)!\,b_2\dotsm b_{k-1}},
	\end{equation}
	the first inequality for $m\geq\Delta_{k-1}$. We substitute~\eqref{eq:refined-known}, with
	$m=n-s-qb_k$, in the terms of~\eqref{eq:refined-rec}.

	\emph{Proof of \textup{(i)}.} Let $n\geq s$, and set $x\coloneqq\frac{n-s+\nabla_{k-1}}{b_k}$,
	so that $m+\nabla_{k-1}=b_k\,(x-q)$. Substituting the upper bound of~\eqref{eq:refined-known}
	in each term of~\eqref{eq:refined-rec}, we get
	\begin{equation*}
		\begin{aligned}
			D_k(n)
			&\leq\frac{b_k^{k-2}}{(k-2)!\,b_2\dotsm b_{k-1}}\sum_{q=0}^{\floor*{\frac{n-s}{b_k}}}(x-q)^{k-2}
			\leq\frac{b_k^{k-2}}{(k-2)!\,b_2\dotsm b_{k-1}}\sum_{q=0}^{\floor{x}}(x-q)^{k-2}\\
			&\leq\frac{b_k^{k-2}\bigl(x+\frac12\bigr)^{k-1}}{(k-1)!\,b_2\dotsm b_{k-1}}
			=\frac{(n-s+\nabla_k)^{k-1}}{(k-1)!\,b_2\dotsm b_k}.
		\end{aligned}
	\end{equation*}
	The second inequality holds because $\frac{n-s}{b_k}\leq x$, so the sum is extended by terms
	with $q\leq x$, which are non-negative. The third one is Lemma~\ref{lem:sum-bounds} with
	$j=k-2$. The final equality is
	$b_k\bigl(x+\frac12\bigr)=n-s+\nabla_{k-1}+\frac{b_k}{2}=n-s+\nabla_k$
	(Lemma~\ref{lem:basic}~(iv)).

	\emph{Proof of \textup{(ii)}.} Let $n\geq s+\Delta_{k-1}$, and set
	$z\coloneqq\frac{n-s-\Delta_{k-1}}{b_k}\geq0$, so that $m-\Delta_{k-1}=b_k\,(z-q)$. For
	$0\leq q\leq\floor{z}$ we have $q\leq z\leq\frac{n-s}{b_k}$, so the term belongs
	to~\eqref{eq:refined-rec}, and $m\geq\Delta_{k-1}$, so the lower bound
	of~\eqref{eq:refined-known} applies to it. Keeping only these terms, and then substituting
	the lower bound in each of them, we get
	\begin{equation*}
		\begin{aligned}
			D_k(n)
			&\geq\sum_{q=0}^{\floor{z}}D_{k-1}(n-s-qb_k)
			>\frac{b_k^{k-2}}{(k-2)!\,b_2\dotsm b_{k-1}}\sum_{q=0}^{\floor{z}}(z-q)^{k-2}\\
			&\geq\frac{b_k^{k-2}\,z^{k-1}}{(k-1)!\,b_2\dotsm b_{k-1}}
			=\frac{(n-s-\Delta_{k-1})^{k-1}}{(k-1)!\,b_2\dotsm b_k}.
		\end{aligned}
	\end{equation*}
	The first inequality holds because the terms that are dropped are non-negative. The second
	one is strict because there is at least the term $q=0$. The third one is the lower bound of
	Lemma~\ref{lem:sum-bounds} with $j=k-2$, without its term $\frac{z^{k-2}}{2}\geq0$. The
	final equality is $b_kz=n-s-\Delta_{k-1}$.
\end{proof}

\begin{rem}\label{rem:refined}
	Proposition~\ref{prop:refined} contains Theorems~\ref{thm:upper} and~\ref{thm:lower}, and is
	sharper. Recall that $0\leq s\leq b_k-a_k$ (Lemma~\ref{lem:recursion}).
	\begin{itemize}
		\item \emph{Upper bound.} If $n\geq s$, then $n-s+\nabla_k\leq n+\nabla_k$, so~(i) implies
		      Theorem~\ref{thm:upper}; if $n<s$, the sum of Lemma~\ref{lem:recursion} is empty and
		      $D_k(n)=0$.
		\item \emph{Lower bound.} By Lemma~\ref{lem:basic}~(iv) (for $k=2$, by $\Delta_2=b_2-a_2$
		      and $\Delta_1=0$),
		      \begin{equation*}
		      	s+\Delta_{k-1}\leq(b_k-a_k)+\Delta_{k-1}=\Delta_k .
		      \end{equation*}
		      Hence every $n\geq\Delta_k$ satisfies $n\geq s+\Delta_{k-1}$, and
		      $n-s-\Delta_{k-1}\geq n-\Delta_k$, so~(ii) implies Theorem~\ref{thm:lower}. Moreover~(ii)
		      gives $D_k(n)>0$ also for the multiples $n$ of $d_k$, if any, with
		      $s+\Delta_{k-1}\leq n<\Delta_k$, where Theorem~\ref{thm:lower} says nothing.
	\end{itemize}
	The refinement concerns only the last step, because any attempt to use it earlier in the
	induction, in order to improve the bounds further, runs into the fact that $s$ depends on
	$n$. Indeed, in the recursion $D_k(n)=\sum_{q}D_{k-1}(n-s-qb_k)$ there is only one integer
	$s$, the one determined by $n$. But to refine also the previous step one has to apply~(i) or~(ii)
	for $D_{k-1}$ to each argument $m=n-s-qb_k$, and each $m$ has its own integer $s(m)$, which
	varies with $q$. Then the terms no longer have the form $(x-q)^{k-2}$ with one and the same
	$x$, and Lemma~\ref{lem:sum-bounds} does not apply. This is why the inductive hypothesis in
	Theorems~\ref{thm:upper} and~\ref{thm:lower} must be uniform in $n$, that is, free of $s$:
	in their proofs $s$ is estimated once and for all by $0\leq s\leq b_{k+1}-a_{k+1}$, and this
	estimate is what produces $\nabla_k$ and $\Delta_k$.

	Finally, the same refinement applies to the polynomial lower bound of
	Theorem~\ref{thm:lower-poly}. For $k\geq3$, in the proof of~(ii) one substitutes in each
	term the bound of Theorem~\ref{thm:lower-poly} for $D_{k-1}$ instead of that of
	Theorem~\ref{thm:lower}. With the polynomials~\eqref{eq:P-def} and the same $z$, the bound
	for each term is then $\frac{P_{k-1}(b_k\,(z-q))}{b_2\dotsm b_{k-1}}$, and from this point on
	the computation is that of Steps~3 and~4 of the proof of Theorem~\ref{thm:lower-poly}, with
	$k-1$ in place of $k$; we do not repeat it. The result is that, for every multiple $n$ of
	$d_k$ with $n\geq s+\Delta_{k-1}$,
	\begin{equation*}
		D_k(n)>\frac{1}{b_2\dotsm b_k}\sum_{i=0}^{k-2}\bbinom{k-2}{i}\frac{(n-s-\Delta_{k-1})^{k-1-i}}{(k-1-i)!},
	\end{equation*}
	that is, Theorem~\ref{thm:lower-poly} with $s+\Delta_{k-1}$ in place of $\Delta_k$; for $k=2$
	this is~(ii).
\end{rem}

\subsection{Solutions of \texorpdfstring{$a_1x_1+\dotsb+a_kx_k\leq n$}{a1x1+...+akxk<=n}}

\begin{defn}\label{dfn:denum-ext}
	For $k\in\N_{\geq1}$, $n\in\N$, let $\widehat{D}_k(n)$ be the number of
	$(x_1,\dotsc,x_k)\in\N^k$ with $a_1x_1+\dotsb+a_kx_k\leq n$.
\end{defn}

To state the bounds for $\widehat{D}_k(n)$ we use the shifted sequence
\begin{equation*}
	a'\coloneqq(a'_1,a'_2,a'_3,\dotsc)=(1,a_1,a_2,\dotsc),
	\qquad\text{that is}\qquad
	a'_1=1,\quad a'_i=a_{i-1}\ (i\geq2),
\end{equation*}
and we denote by $d'_i$, $b'_i$, $\Delta'_k$, $\nabla'_k$ and $\bbinom[b']{m}{i}$ the objects of
Definitions~\ref{dfn:denumerant} and~\ref{dfn:constants} and Section~\ref{s:prelim} built from $a'$ instead of $a$.
Since $a'_1=1$, Lemma~\ref{lem:basic}~(v), applied to $a'$, gives $d'_i=1$, $b'_i=a'_i$ for
every $i$, and $\Delta'_k=0$; hence
\begin{equation*}
	b'_2\dotsm b'_{k+1}=a_1\dotsm a_k,
	\qquad
	\nabla'_{k+1}=b'_2+\frac12\sum_{i=3}^{k+1}b'_i=a_1+\frac12\sum_{i=2}^{k}a_i,
	\qquad
	\bbinom[b']{k-1}{i}=e_i\Bigl(\frac{a_2}{2},\dotsc,\frac{a_k}{2}\Bigr).
\end{equation*}

\begin{cor}\label{cor:denum-ext}
	Let $k\in\N_{\geq1}$ and $n\in\N$, and set $V\coloneqq d_k\floor*{\frac{n}{d_k}}$. Then
	\begin{equation*}
		\frac{V^k}{k!\,a_1\dotsm a_k}
		\leq\frac{1}{a_1\dotsm a_k}\sum_{i=0}^{k-1}\bbinom[b']{k-1}{i}\frac{V^{k-i}}{(k-i)!}
		<\widehat{D}_k(n)
		\leq\frac{(V+\nabla'_{k+1})^k}{k!\,a_1\dotsm a_k}.
	\end{equation*}
\end{cor}
\begin{proof}
	\emph{Reduction to $V$.} By Lemma~\ref{lem:basic}~(vi), every value $a_1x_1+\dotsb+a_kx_k$ is a
	multiple of $d_k$, and $V$ is the largest multiple of $d_k$ not exceeding $n$. Hence
	$a_1x_1+\dotsb+a_kx_k\leq n$ if and only if $a_1x_1+\dotsb+a_kx_k\leq V$, that is
	$\widehat{D}_k(n)=\widehat{D}_k(V)$.

	\emph{A slack variable.} The map
	\begin{equation*}
		(x_1,\dotsc,x_k)\longmapsto\bigl(V-a_1x_1-\dotsb-a_kx_k,\ x_1,\dotsc,x_k\bigr)
	\end{equation*}
	is a bijection between the solutions in $\N^k$ of $a_1x_1+\dotsb+a_kx_k\leq V$ and the
	solutions in $\N^{k+1}$ of $x_0+a_1x_1+\dotsb+a_kx_k=V$, its inverse being
	$(x_0,x_1,\dotsc,x_k)\mapsto(x_1,\dotsc,x_k)$. The latter equation is
	$a'_1x_0+a'_2x_1+\dotsb+a'_{k+1}x_k=V$, so $\widehat{D}_k(V)$ is the denumerant
	$D_{k+1}(V)$ of the sequence $a'$.

	\emph{Conclusion.} Since $k+1\geq2$, $d'_{k+1}=1$ and $V\geq0=\Delta'_{k+1}$,
	Theorems~\ref{thm:upper} and~\ref{thm:lower-poly}, applied to the sequence $a'$ with $k+1$
	in place of $k$, give
	\begin{equation*}
		\frac{1}{b'_2\dotsm b'_{k+1}}\sum_{i=0}^{k-1}\bbinom[b']{k-1}{i}\frac{(V-\Delta'_{k+1})^{k-i}}{(k-i)!}
		<\widehat{D}_k(V)
		\leq\frac{(V+\nabla'_{k+1})^{k}}{k!\,b'_2\dotsm b'_{k+1}},
	\end{equation*}
	which, with $\Delta'_{k+1}=0$ and $b'_2\dotsm b'_{k+1}=a_1\dotsm a_k$, are the two inequalities
	on the right. The inequality on the left holds because, by Lemma~\ref{lem:sym-poly}, applied to
	$a'$, the term $i=0$ of the sum is $\frac{V^k}{k!}$ and the other terms are non-negative.
\end{proof}

\begin{rem}\label{rem:denum-ext}
	The bounds of Corollary~\ref{cor:denum-ext} depend on $n$ only through $V$, as
	$\widehat{D}_k(n)$ does.
\end{rem}

\bigskip
\noindent\textbf{Acknowledgements.} The author used an AI assistant (Claude, Anthropic) to revise
the English wording and the organization of the text, and to check the formulas and the
numerical examples. All results and proofs are due to the author, who takes full responsibility
for the content of the paper.

\bibliographystyle{aomalpha}
\bibliography{denumerant_bounds_I}

\end{document}